\documentclass[11pt,a4paper,3p]{elsarticle}
\usepackage[utf8]{inputenc}
\usepackage[T1]{fontenc}
\usepackage[english]{babel}
\usepackage{enumerate}
\usepackage{graphicx} 
\usepackage{booktabs}
\usepackage[lofdepth,lotdepth]{subfig}
\usepackage{mathtools}
\usepackage{amsmath,bm}
\usepackage{amssymb}
\usepackage{fancyhdr}
\usepackage{fourier} 
\usepackage{pgfplots} 
\usepackage{cancel}
\usepackage{graphicx}        
\newtheorem{remark}{Remark}[section]

\newcommand{\vect}[1]{\mathbf{#1}}
\newcommand{\tdomain}{\Omega_1(t) \cup \Omega_2(t)}

\newcommand{\bfx}{\mathbf{x}}
\newcommand{\bfv}{\mathbf{v}}
\newcommand{\bfu}{\mathbf{u}}
\newcommand{\bfw}{\mathbf{w}}

\newcommand{\bfn}{\mathbf{n}}
\newcommand{\bfg}{\mathbf{g}}
\newcommand{\bff}{\mathbf{f}}
\newcommand{\bfI}{\mathbf{I}}
\newcommand{\bfeps}{\bm{\epsilon}}
\newcommand{\mcK}{\mathcal{K}}
\newcommand{\mcN}{\mathcal{N}}
\newcommand{\lsf}{\phi}
\newcommand{\rf}{\hat{\phi}}
\newcommand{\vel}{\mathbf{\beta}}

\newcommand{\R}{\mathbb{R}}

\newcommand{\Po}{\mathbb{P}}

\newcommand{\enstq}[2]{\left\{#1\mathrel{}\middle|\mathrel{}#2\right\}}

\newcommand{\norm}[1]{\left\Vert #1\right\Vert}
\newcommand{\sign}[1]{\text{sign}\left( #1 \right)}
\renewcommand{\div}[1]{\nabla \cdot #1}

\newcommand{\grad}[1]{\nabla #1}

\newcommand{\diff}{\mathop{ }\mathopen{ }\mathrm{d}}
\newcommand{\restreinta}[1]{\mathclose{}|\mathopen{}_{#1}}
\newcommand{\prodscal}[2]{\left( #1\ ,\ #2\right)}
\newcommand{\jump}[1]{\llbracket #1 \rrbracket}
\newcommand{\dOmega}{\partial \Omega}
\newcommand{\dw}{\partial \Omega}

\newcommand{\w}{k} 
\newcommand{\mcF}{\mathcal{F}}
\newcommand{\mapp}{\theta}
\newcommand{\meancurv}{\mathcal{H}}
\begin{document}

\begin{frontmatter}

\title{A cut finite element method for incompressible two-phase Navier--Stokes flows }

\author[1]{Thomas Frachon%
}
\ead{frachon@kth.se}
\address[1]{Department of Mathematics, KTH Royal Institute of Technology,
SE-100 44 Stockholm, Sweden}

\author[1]{Sara Zahedi\corref{cor1}%
}
\ead{sara.zahedi@math.kth.se}

\cortext[cor1]{Corresponding author}

\begin{abstract}
We present a space-time Cut Finite Element Method (CutFEM) for the time-dependent Navier--Stokes equations involving two immiscible incompressible fluids with different viscosities, densities, and with surface tension. 
The numerical method is able to accurately capture the strong discontinuity in the pressure and the weak discontinuity in the velocity field across evolving interfaces without re-meshing processes or regularization of the problem. We combine the strategy proposed in [P. Hansbo, M. G. Larson, S. Zahedi, Appl. Numer. Math. 85 (2014), 90--114] for the Stokes equations with a stationary interface and the space-time strategy presented in [P. Hansbo, M. G. Larson, S. Zahedi, Comput. Methods Appl. Mech. Engrg. 307 (2016), 96--116].  
We also propose a strategy for computing high order approximations of the surface tension force by computing a stabilized mean curvature vector. The presented space-time CutFEM uses a fixed mesh but includes stabilization terms that control the condition number of the resulting system matrix independently of the position of the interface, ensure stability and a convenient implementation of the space-time method based on quadrature in time.  Numerical experiments in two and three space dimensions show that the numerical method is able to accurately capture the discontinuities in the pressure and the velocity field across evolving interfaces without requiring the mesh to be conformed to the interface and with good stability properties.
\end{abstract}


\end{frontmatter}

\section{Introduction}
Today computer simulation provides valuable insights of two-phase flow phenomena and is an important tool in studies of such flow problems. For reliable simulations, accurate and robust computational techniques are needed and much effort is directed to their development, see e.g. \cite{UnTr92, SuSmOs94, FrBe10, GrRe11, GaTo12, BGN15, GiHyFe18} and references therein.

From a computational point of view one of the main challenges in two-phase flow simulations is that Partial Differential Equations (PDEs) need to be solved in evolving domains. These evolving domains are defined by interfaces that separate the immiscible fluids.
When an interface undergoes large deformations, for example when topological changes such as drop-breakup or coalescence occur, the re-meshing and interpolation that is required by standard Finite Element Methods (FEM), as well as standard finite difference schemes, becomes both cumbersome and expensive, especially in three space dimensions. Therefore, for simulations of problems where the interface may undergo large deformations so called fixed-grid flow solvers are desirable. These flow solvers must be able to accurately approximate discontinuities in the solution across deforming interfaces without conforming the mesh to these interfaces. Several strategies exist, see e.g. \cite{Win07, GrRe07, FrBe10, ScRaGrWa15, KrBoSiVo}.

The Cut Finite Element Method (CutFEM) is a robust and accurate unfitted finite element method which, contrary to standard FEM, allows the evolving geometry to be arbitrarily located with respect to a fixed background mesh, but has the same order of accuracy and scaling of the condition number with respect to the mesh size as standard FEM. Discontinuities in the solution across an interface are accurately captured by building the solution from two solutions, one on each subdomain separated by the interface, and then glue the solutions at the interface by weakly imposing the interface conditions in the variational formulation, see e.g. \cite{HaHa02, BBH09, HaLaZa14, BuClHaLaMa15}. Stabilization terms that are weak enough to not destroy the optimal convergence order but strong enough to ensure well-posedness of the resulting algebraic system of equations independently of how the interface cuts the background mesh, are added to the variational formulation \cite{Bur10, BH12, WZKB}. Stabilization terms may also be added by other reasons for example to improve the accuracy in the computation of the mean curvature vector as in \cite{HanLarZah15}, or to obtain a stable discretization \cite{HaLaZa14, HLZ16, Zah18}, see also  Section \ref{sec:weakform_nummeth} and \ref{sec:surften}. 
In CutFEM the discretization of the PDE is independent of the numerical representation of the interface, and different techniques for representing and evolving the interface can be used.

In this work we consider the Navier--Stokes equations governing the fluid motion of two immiscible incompressible fluids. For this problem we propose a second order accurate space-time CutFEM which accurately captures the discontinuities, in both parameters and in the solution, across the evolving interface. Our method avoids both re-meshing processes and regularization of discontinuities. It is built on the CutFEM we proposed and analyzed in \cite{HaLaZa14} for the Stokes equations involving two stationary immiscible incompressible fluids and the space-time CutFEM we proposed in \cite{HLZ16} for convection-diffusion equations in time dependent domains. In the proposed space-time method, we use a CutFEM in space based on inf-sup stable elements and we use discontinuous piecewise polynomials in time, stabilization terms are added that ensure good stability properties and a convenient implementation using quadrature rules in time to directly approximate space-time integrals in the variational formulation.  
We also propose a method for computing high order accurate approximations of the surface tension force by stabilizing the $L^2$ projection and computing a stabilized mean curvature vector approximation. Our method is an extension of the method we presented in \cite{HanLarZah15} for piecewise linear interface approximations. In this work we have used a level set method to represent the interface and the mapping proposed in \cite{Le16} to transform integrals on high order approximations of the interface, implicitly defined through the level set function, to integrals on a piecewise linear approximation of the interface. We emphasize that, provided a numerical method for representing and evolving the interface exist, the proposed space-time CutFEM has a convenient implementation as it is built on a stationary implementation of CutFEM at discrete time instances.

The rest of the paper is outlined as follows. In Section 2, we state the governing equations and propose a variational formulation in which the physical interface conditions are imposed weakly. The numerical method, a space-time CutFEM based on quadrature rules in time, is presented in Section 3. We describe the level set method we have used in Section 4 and in Section 5 we present a high order method for computing the surface tension force. Numerical examples are shown in Section 6. We summarize this work in Section 7.

\section{Governing equations}
\begin{figure}
\begin{center}
\includegraphics[width=0.35\textwidth]{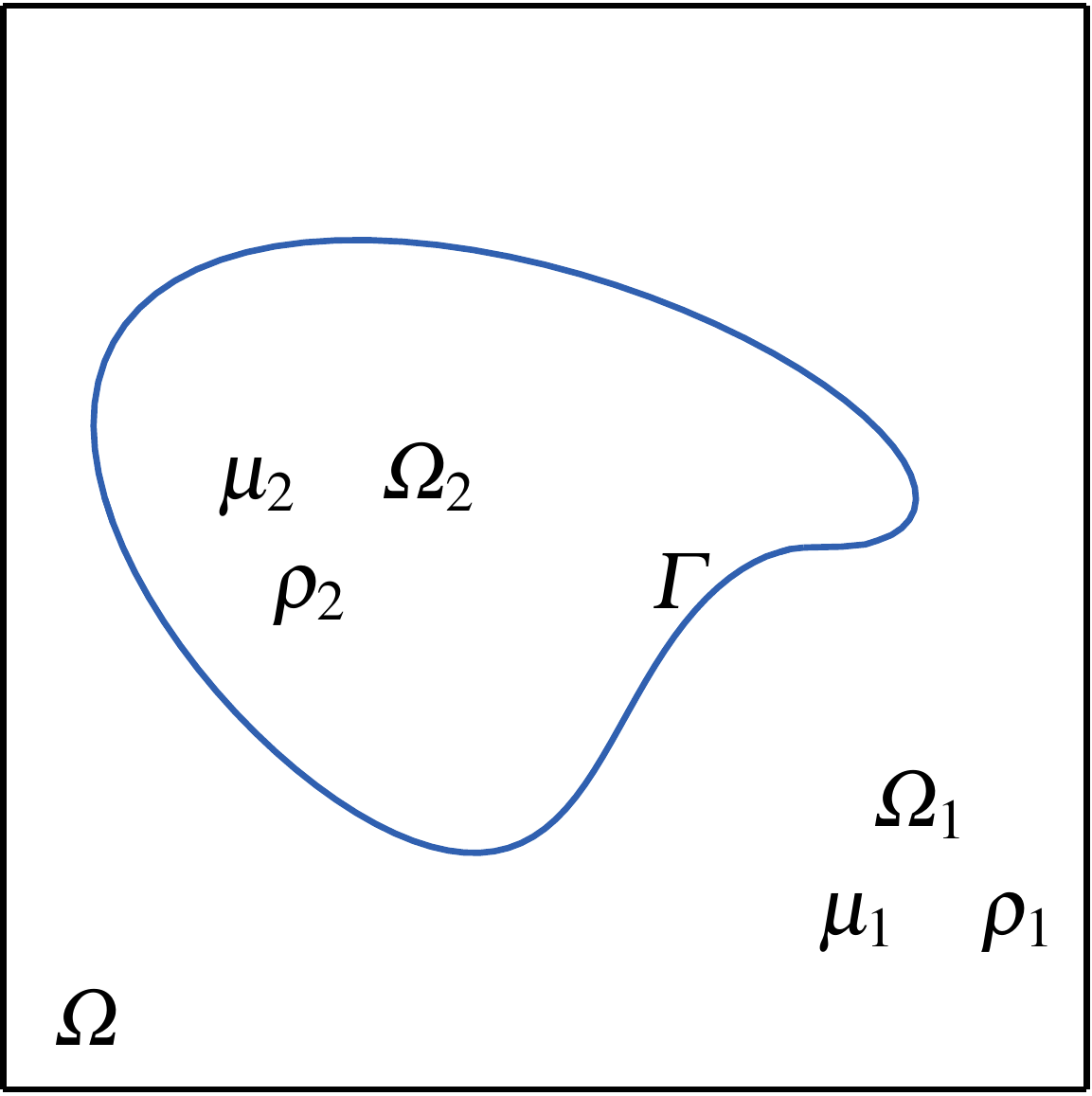} 
\caption{Illustration of the domain $\Omega \in \R^2$ and the two subdomains $\Omega_i(t)$, $i=1,2$ occupied by immiscible fluids with different material properties $\rho_i$ and $\mu_i$, separated by an interface $\Gamma$.\label{fig:domain}}
\end{center}
\end{figure}

We consider the dynamics of two immiscible incompressible fluids with different material properties contained in a bounded domain $\Omega$ in $\R^d$, $d=2,3$, with a convex polygonal boundary $\dw$. During time $t$ in a time interval $I= [0,T]$, the two fluids occupy time dependent subdomains $\Omega_i(t) \subset \Omega, i=1,2$, such that $\bar{\Omega} = \bar{\Omega}_1(t) \cup \bar{\Omega}_2(t)$ and $\Omega_1(t) \cap \Omega_2(t) = \emptyset$. Furthermore, the two immiscible fluids are separated by a sufficiently smooth interface defined by $\Gamma(t) = \dw_1(t) \cap \dw_2(t)$. We assume that $\Omega_2(t)$ is the domain enclosed by $\Gamma(t)$. See Fig.~\ref{fig:domain} for an illustration in two dimensions ($d=2$). We assume the dynamics of the fluids is governed by the incompressible Navier--Stokes equations. 

\subsection{The two-phase incompressible Navier--Stokes equations}
For $t\in I$ the motion of the two immiscible incompressible fluids with viscosities $\mu_i$, densities $\rho_i$, $i=1,2$, and surface tension is assumed to be governed by the following standard model: 
\begin{alignat}{2}
\rho  \left( \partial_t \bfu + (\bfu \cdot \grad){\bfu} \right)  
	- \div{\left( 2\mu \bfeps(\bfu) - p \bfI \right)}  &= \bff  & &\quad \text{in} \ \tdomain,  \label{eq : Navier-Stokes 1} \\
\div{\bfu} & =  0   & &\quad \text{in} \ \tdomain, \label{eq : Incompressibility}\\
\jump{\bfu} & =  0 & &\quad \text{on} \ \Gamma(t), \label{eq : Continuity velocity}\\
\jump{(2\mu \bfeps(\bfu) - p \mathbf{I}) \vect{n}} & =  \sigma \kappa \vect{n} & & \quad \text{on} \ \Gamma(t), \label{eq : Surface tension 1}\\
\bfu(t,\bfx) & =  \bfg(t,\bfx) &  &\quad \text{on} \ \partial \Omega, \label{eq : Boundary condition}\\
\bfu(0,\bfx) & =  \bfu^0(\bfx) & &\quad \text{in} \ \Omega_1(0) \cup \Omega_2(0).\label{eq : Initial condition}
\end{alignat}
Here, $\bfu : I \times \Omega \rightarrow \R^d$ is the velocity field, $p : I \times \Omega \rightarrow \R$ is the pressure,  $\partial_t = \frac{\partial}{\partial t}$, $\bfeps(\bfu) = (\grad{\bfu} + \grad{\bfu}^T)/2$ is the strain rate tensor, $\mu = \mu_i > 0$ in $\Omega_i(t)$, $\rho = \rho_i > 0$, in $\Omega_i(t)$, $i=1,2$, $\sigma$ is the surface tension coefficient, $\kappa$ is the mean curvature of the interface, $\bfn$ is the outward directed (with respect to $\Omega_1$) unit normal to $\partial \Omega_1$, $\jump{v} = (v_1 - v_2)\restreinta{\Gamma}$ is the jump, where $v_i = v \restreinta{\Omega_i}, i=1,2$, $\bff: I \times \Omega \rightarrow \R^d$ is a given external volume force (e.g. the gravitational force), $\bfg$ is a given function such that $\int_{\partial \Omega} \bfg \cdot \bfn \diff s=0$, and $\bfu^0:\Omega \rightarrow \R^d$ is given and defines the initial condition. The initial configuration of the interface $\Gamma(0)$ and thus the subdomains $\Omega_i(0)$, $i=1,2$ is also given. For simplicity, we will assume that the interface does not intersect the boundary of the domain $\Omega$ during time $t\in I$. 
Note that in this formulation, if $(\bfu, p)$ is a solution then $(\bfu,p+c)$, $c \in \R$ is also a solution. 

We may have other type of boundary conditions, for example the mixed boundary conditions
\begin{align}\label{eq:mixedBC}
\bfu \cdot \bfn&=\bfg \cdot \bfn, \nonumber \\
(2\mu \bfeps(\bfu)\vect{n})\cdot \vect{\tau}&=0, \quad  \vect{\tau} \cdot \bfn=0,
\end{align}
on $\partial \Omega$ or parts of $\partial \Omega$. 

In the following subsection we formulate a variational formulation of problem \eqref{eq : Navier-Stokes 1}-\eqref{eq : Initial condition} in which the interface and the boundary conditions are imposed weakly.

\subsection{A weak formulation}

Recall the Sobolev spaces
\begin{equation}
L^2(U)=\enstq{v}{\int_{U} |v|^2 \diff \bfx=\| v \|_{L^2(U)}^2 < \infty},
\end{equation}
and

\begin{equation}
H^1(U)=\enstq{v}{v \in L^2(U), \grad v \in L^2(U)},
\end{equation}
where $U$ is a domain in $\R^d$. We will use the notation $\prodscal{v}{w}_U=\int_U v(\bfx) w(\bfx) \diff \bfx$ for the $L^2$ inner product on $U$ (similarly for inner products in $[L^2(U)]^d$) and 
\begin{equation}
\prodscal{v}{w}_{U_1 \cup U_2}=\sum_{i=1}^2\prodscal{v}{w}_{U_i}.
\end{equation}

Now introduce the spaces 
\begin{align*}
V = & \left[H^1(\tdomain) \right]^d, \\
Q =& \enstq{q \in L^2(\Omega) }{\prodscal{\mu^{-1}q}{1}_{\tdomain} = 0}.
\end{align*}
For $t\in I$  we formulate the following weak formulation of problem \eqref{eq : Navier-Stokes 1}-\eqref{eq : Initial condition}: find $(\bfu,p)\in V \times Q$ such that 
\begin{equation}\label{eq:weakformNS}
\prodscal{\rho \partial_t \bfu}{\bfv}_{\tdomain } + \prodscal{\rho(\bfu \cdot \grad){\bfu}}{\bfv}_{\tdomain} +a(t, \bfu,\bfv)-b(t,\bfv,p)+b(t,\bfu,q)=l(t,v,q) 
\end{equation}
for all $(\bfv,q) \in V \times Q$ and $\bfu(0,\bfx)  =  \bfu^0(\bfx)$ in $\Omega_1(0) \cup \Omega_2(0)$.
Here
\begin{align}\label{eq:forma}
	a(t,\bfu,\bfv) & = 
        \prodscal{2\mu\varepsilon(\bfu)}{\varepsilon(\bfv)}_{\tdomain} \nonumber \\
&	- \prodscal{\lbrace 2\mu \varepsilon(\bfu)  \bfn \rbrace}{ \jump{\bfv}}_{\Gamma(t)}
	-  \prodscal{ \jump{\bfu}}{\lbrace 2\mu \varepsilon(\bfv)  \bfn \rbrace}_{\Gamma(t)}
	+ \prodscal{\lambda_{\Gamma} \jump{\bfu}}{\jump{\bfv}}_{\Gamma(t)} \nonumber\\
&	- \prodscal{2\mu \varepsilon(\bfu)  \bfn}{\bfv}_{\dOmega}
	- \prodscal{\bfu}{2\mu \varepsilon(\bfv)  \bfn}_{\dOmega}	
	+ \prodscal{\lambda_{\dOmega}  \bfu}{\bfv}_{\dOmega},
\end{align}
$b(t,\bfv,q)=b^1(t,\bfv,q)$ or $b(t,\bfv,q)=b^2(t,\bfv,q)$  with
\begin{align}
	b^1(t,\bfv,q) & = \prodscal{\div \bfv}{q}_{\tdomain}  
	- \prodscal{\jump{\bfv \cdot \bfn}}{ \{ q \}}_{\Gamma(t)}	
	- \prodscal{\bfv \cdot \bfn}{q}_{\partial \Omega}, 
\label{eq:formb1}\\
 	b^2(t,\bfv,q) & = 
	-  \prodscal{\bfv}{\grad{q}}_{\tdomain} 
	+\prodscal{\jump{q}}{ \langle \bfv \cdot \bfn \rangle}_{\Gamma(t)},\label{eq:formb2}
\end{align}
both forms are mathematically equivalent,
\begin{align}\label{eq:forml}
	l(t,\bfv,q) & = \prodscal{\mathbf{f}}{\bfv}_{\tdomain} 
	+ \prodscal{\sigma \kappa \bfn}{\langle \bfv \rangle}_{\Gamma(t)} 
	- \prodscal{\bfg }{ 2\mu \varepsilon(\bfv)  \bfn}_{\dOmega} 
        +\prodscal{\lambda_{\dOmega} \vect{g} }{ \bfv}_{\dOmega} 
        - \prodscal{\bfg \cdot \bfn }{q}_{\dOmega}, 
\end{align}
and 
\begin{equation}\label{eq:averagop}
\{f \}=\w_1f_1+\w_2f_2, \qquad  \langle f \rangle=\w_2f_1+\w_1f_2,
\end{equation}
where the weights $\w_1$ and $\w_2$ are real numbers satisfying $\w_1+\w_2=1$ and $f_i = f \restreinta{\Omega_i}$.

We now derive the given variational formulation using a variant of Nitsche's method \cite{Nit}. For $t\in I$ assume that $(\bfu,p)$ with $p$ such that $\prodscal{\mu^{-1}p}{1}_{\tdomain} = 0$ is a sufficiently smooth solution to \eqref{eq : Navier-Stokes 1}-\eqref{eq : Initial condition} (other conditions than $\prodscal{\mu^{-1}q}{1}_{\tdomain} = 0$ can also be used to fix the constant in the pressure).   
Multiply both sides of equation \eqref{eq : Navier-Stokes 1}  with a smooth vector field $\bfv$, integrate in each subdomain $\Omega_i(t)$, $i=1,2$ and apply integration by parts in the different subdomains, to arrive at 
\begin{align} \label{eq:partialinteq1}
&	\prodscal{\rho \partial_t \bfu}{\bfv}_{\tdomain} 
+	\prodscal{\rho (\bfu \cdot \grad){\bfu}}{\bfv}_{\tdomain} 
+	\prodscal{2\mu\varepsilon(\bfu)}{\varepsilon(\bfv)}_{\tdomain}  
-	\prodscal{2\mu \varepsilon(\bfu)\bfn}{\bfv}_{\dOmega}   \nonumber \\
-&	\prodscal{p}{\div{\bfv}}_{\tdomain} 
+	\prodscal{p}{\bfn\cdot \bfv}_{\dOmega} 
-       \int_{\Gamma(t)} \jump{((2\mu \varepsilon(\bfu) - p \bfI)\bfn)\bfv} \diff s 
=	\prodscal{f}{\bfv}_{\tdomain}.
\end{align}

We used the following partial integration rules
\begin{align}
\prodscal{-\div (2\mu\varepsilon(\bfu))}{\bfv}_{\Omega_i}=
\prodscal{2\mu\varepsilon(\bfu)}{\varepsilon(\bfv)}_{\Omega_i}
- \prodscal{2\mu\varepsilon(\bfu)\bfn_i}{\bfv}_{\partial \Omega_i}, 
\label{eq:parintrul2}\\
\prodscal{\nabla p}{\bfv}_{\Omega_i}=-\prodscal{p}{\div {\bfv}}_{\Omega_i} +\prodscal{p}{\bfn_i \cdot \bfv}_{\partial \Omega_i}, \label{eq:parintrul3}
\end{align}
where $\bfn_i$ is the outward directed (with respect to $\Omega_i$) unit normal to $\partial \Omega_i$. Note that with our definition for $\bfn$ we have $\bfn=\bfn_1$ and $\bfn_2=-\bfn$. 

It can easily be checked that 
\begin{equation}\label{eq:relationaverag}
\jump {f g}=\{ f\} \jump{g}+ \jump{f} \langle g \rangle,
\end{equation} 
holds for the averaging operators in equation \eqref{eq:averagop} since $\w_1+\w_2=1$. Using equation \eqref{eq:relationaverag} (with $f=(2\mu \varepsilon(\bfu) -p \mathbf{I}) \bfn$ and $g=\bfv$), we can rewrite the last integral on the left hand side of equation \eqref{eq:partialinteq1} in terms of the jump in normal stress and apply the interface condition \eqref{eq : Surface tension 1}: 
\begin{equation}\label{eq:rewriteinterfT}
\int_{\Gamma}  \jump{((2\mu \varepsilon(\bfu) -p \mathbf{I}) \bfn) \bfv} \diff s 
= \int_{\Gamma}  \{(2\mu \varepsilon(\bfu) -p \mathbf{I}) \bfn \} \jump{\bfv} \diff s 
+ \int_{\Gamma} \underbrace{ \jump{(2\mu \varepsilon(\bfu) -p \mathbf{I}) \bfn}}_{\sigma \kappa \bfn} \langle \bfv \rangle \diff s.
\end{equation}

Now multiplying equation \eqref{eq : Incompressibility} with a smooth function $q$ and integrating in the different subdomains we have
\begin{align} \label{eq:formbeforeinterfaceterm}
&	\prodscal{\rho \partial_t \bfu}{\bfv}_{\tdomain} 
+	\prodscal{\rho (\bfu \cdot \grad){\bfu}}{\bfv}_{\tdomain} 
+	\prodscal{2\mu\varepsilon(\bfu)}{\varepsilon(\bfv)}_{\tdomain}  
-	\prodscal{2\mu \varepsilon(\bfu)\bfn}{\bfv}_{\dOmega} \nonumber \\
-&	\prodscal{p}{\div {\bfv}}_{\tdomain} 
+	\prodscal{p}{\bfn \cdot \bfv}_{\dOmega} 
-	\prodscal{\lbrace (2\mu \varepsilon(\bfu) -pI)\bfn \rbrace}{ \jump{\bfv}}_{\Gamma(t)}    +(\div \bfu,q)_{\tdomain}
= \nonumber\\
=&	 \prodscal{f}{\bfv}_{\tdomain} + \prodscal{\sigma \kappa \bfn}{\langle \bfv \rangle}_{\Gamma(t)} .
\end{align}

From the interface condition \eqref{eq : Continuity velocity} and the boundary condition \eqref{eq : Boundary condition}, it follows that
\begin{align}
-\prodscal{ \jump{\bfu}}{\lbrace 2\mu \varepsilon(\bfv)\bfn \rbrace}_{\Gamma(t)}  =0, \label{eq:termadd1}\\
-\prodscal{\bfu}{2\mu \varepsilon(\bfv) \bfn}_{\dOmega}
+\prodscal{\bfg}{2\mu \varepsilon(\bfv) \bfn}_{\dOmega}=0, \label{eq:termadd2}\\
\prodscal{ \jump{\bfu}}{\lbrace -qI\bfn \rbrace}_{\Gamma(t)}  =0, \\
\prodscal{\bfu}{(-q\bfI)\bfn}_{\dOmega}
-\prodscal{\bfg}{(-q\bfI)\bfn}_{\dOmega}=0, \\
\prodscal{\lambda_\Gamma\jump{\bfu}}{\jump{\bfv}}_{\Gamma(t)} =0, \\
\prodscal{\lambda_{\partial \Omega}\bfu}{\bfv}_{\partial \Omega}-\prodscal{\lambda_{\partial \Omega}\bfg}{\bfv}_{\partial \Omega} =0,
\end{align}
for $\lambda_\Gamma \in L^{\infty}(\Gamma)$ and $\lambda_{\partial \Omega} \in L^{\infty}(\partial \Omega)$.  The last two terms are the Nitsche penalty terms. Adding these expressions into \eqref{eq:formbeforeinterfaceterm} yields the given weak formulation with $b=b^1$. Starting from $b^1$, integrating by parts on each subdomain using equation \eqref{eq:parintrul3}, and equation \eqref{eq:relationaverag} one obtains $b^2$. 

\begin{remark}
The signs in equation \eqref{eq:termadd1} and \eqref{eq:termadd2} were chosen to get a symmetric form $a(t,\bfu,\bfv)$ if the signs are changed we instead obtain the nonsymmetric form  
\begin{align}\label{eq:formanonsym}
	a(t,\bfu,\bfv) & = 
        \prodscal{2\mu\varepsilon(\bfu)}{\varepsilon(\bfv)}_{\tdomain} \\
&	- \prodscal{\lbrace 2\mu \varepsilon(\bfu) \cdot \bfn \rbrace}{ \jump{\bfv}}_{\Gamma(t)}
	+  \prodscal{ \jump{\bfu}}{\lbrace 2\mu \varepsilon(\bfv) \cdot \bfn \rbrace}_{\Gamma(t)}
	+ \prodscal{\lambda_{\Gamma(t)} \jump{\bfu}}{\jump{\bfv}}_{\Gamma(t)}\\
&	- \prodscal{2\mu \varepsilon(\bfu) \cdot \bfn}{\bfv}_{\dOmega}
	+ \prodscal{\bfu}{2\mu \varepsilon(\bfv) \cdot \bfn}_{\dOmega}	
	+ \prodscal{\lambda_{\dOmega}  \bfu}{\bfv}_{\dOmega}
\end{align}
and
\begin{align}\label{eq:formlnonsym}
	l(t,\bfv,q) & = \prodscal{\mathbf{f}}{\bfv}_{\tdomain} 
	+ \prodscal{\sigma \kappa \bfn}{\langle \bfv \rangle}_{\Gamma(t)} 
	+\prodscal{\vect{g} }{ 2\mu \varepsilon(\bfv) \bfn}_{\dOmega} 
        +\prodscal{\lambda_{\dOmega} \vect{g} }{ \bfv}_{\dOmega} 
        - \prodscal{\vect{g} \cdot \bfn }{q}_{\dOmega}. 
\end{align}
In \cite{BoiBur16} a finite element method based on the nonsymmetric version of Nitsche's method is studied for linear elasticity and it is shown that the penalty term can be eliminated. However, as a consequence of the nonsymmetric method the convergence order in the $L^2$-error is often suboptimal, see \cite{BoiBur16} and references therein.
\end{remark}

\begin{remark}
In case the stationary Stokes equations are considered with a stationary interface the first two terms in equation \eqref{eq:weakformNS} vanish, the forms $a$, $b$, and $l$ do not depend on time, and all integrals are evaluated on stationary subdomains $\Omega_i$ and $\Gamma$. The weak formulation reduces to the same formulation as proposed in \cite{HaLaZa14}. 
\end{remark}

\begin{remark}
We emphasize that the second term in equation \eqref{eq:forml} contains the jump in normal stress across the interface, see equation \eqref{eq:rewriteinterfT}. Thus, if the interface condition \eqref{eq : Surface tension 1} changes, the term $\prodscal{\sigma \kappa \bfn}{\langle \bfv \rangle}_{\Gamma(t)}$ should be modified. This is the case if surfactants are present. In that case the surface tension coefficient $\sigma$ is not constant and the interface condition \eqref{eq : Surface tension 1} changes to   
\begin{equation}
 \jump{(2\mu \bfeps(\bfu) - p \mathbf{I}) \vect{n}}  =  \sigma \kappa \vect{n} -\nabla_\Gamma\sigma. 
\end{equation}
Hence the second term in equation \eqref{eq:forml} would be $\prodscal{\sigma \kappa \bfn-\nabla_\Gamma\sigma}{\langle \bfv \rangle}_{\Gamma(t)}$.
\end{remark}

\begin{remark}\label{rem:mixedBC}
In case the boundary $\partial \Omega$ is split into subsets $\partial \Omega_D$ and $\partial \Omega_M$ where on $\partial \Omega_D$ we have Dirichlet boundary conditions on $\bfu$ and on $\partial \Omega_M$ we have the mixed boundary conditions in equation \eqref{eq:mixedBC}, then the forms $a(t,\bfu,\bfv)$ and $l(t,\bfv,q)$ change to
\begin{align}\label{eq:formamixedBC}
	a(t,\bfu,\bfv) & = 
        \prodscal{2\mu\varepsilon(\bfu)}{\varepsilon(\bfv)}_{\tdomain} \nonumber \\
&	- \prodscal{\lbrace 2\mu \varepsilon(\bfu)  \bfn \rbrace}{ \jump{\bfv}}_{\Gamma(t)}
	-  \prodscal{ \jump{\bfu}}{\lbrace 2\mu \varepsilon(\bfv)  \bfn \rbrace}_{\Gamma(t)}
	+ \prodscal{\lambda_{\Gamma(t)} \jump{\bfu}}{\jump{\bfv}}_{\Gamma(t)} \nonumber\\
&	- \prodscal{2\mu \varepsilon(\bfu)  \bfn}{\bfv}_{\dOmega_D}
	- \prodscal{\bfu}{2\mu \varepsilon(\bfv)  \bfn}_{\dOmega_D}	
	+ \prodscal{\lambda_{\dOmega}  \bfu}{\bfv}_{\dOmega_D}
\nonumber \\
&        - \prodscal{2\mu \varepsilon(\bfu)  \bfn}{\bfv \cdot \bfn}_{\dOmega_M}
	- \prodscal{\bfu \cdot \bfn}{2\mu \varepsilon(\bfv)  \bfn}_{\dOmega_M}	
	+ \prodscal{\lambda_{\dOmega}  \bfu \cdot \bfn}{\bfv\cdot \bfn}_{\dOmega_M},
\end{align}
and
\begin{align}\label{eq:formldiffBC}
	l(t,\bfv,q) & = \prodscal{\mathbf{f}}{\bfv}_{\tdomain} 
	+ \prodscal{\sigma \kappa \bfn}{\langle \bfv \rangle}_{\Gamma(t)} 
	- \prodscal{\vect{g} }{ 2\mu \varepsilon(\bfv)  \bfn}_{\dOmega_D} 
        +\prodscal{\lambda_{\dOmega} \vect{g} }{ \bfv}_{\dOmega_D} 
        - \prodscal{\vect{g} \cdot \bfn }{q}_{\dOmega_D \cup \dOmega_M}
\nonumber \\
&  - \prodscal{\vect{g} \cdot \bfn }{ 2\mu \varepsilon(\bfv)  \bfn}_{\dOmega_M} 
        +\prodscal{\lambda_{\dOmega} \vect{g} \cdot \bfn }{ \bfv \cdot \bfn}_{\dOmega_M}. 
\end{align}
\end{remark}

\section{A space-time cut finite element method based on quadrature in time}\label{sec:nummeth}
We will present a numerical method for the two-phase incompressible Navier--Stokes equations based on the weak formulation we derived in the previous section.  

\subsection{Mesh}\label{sec:mesh}
\begin{figure}\centering
\includegraphics[scale=0.45]{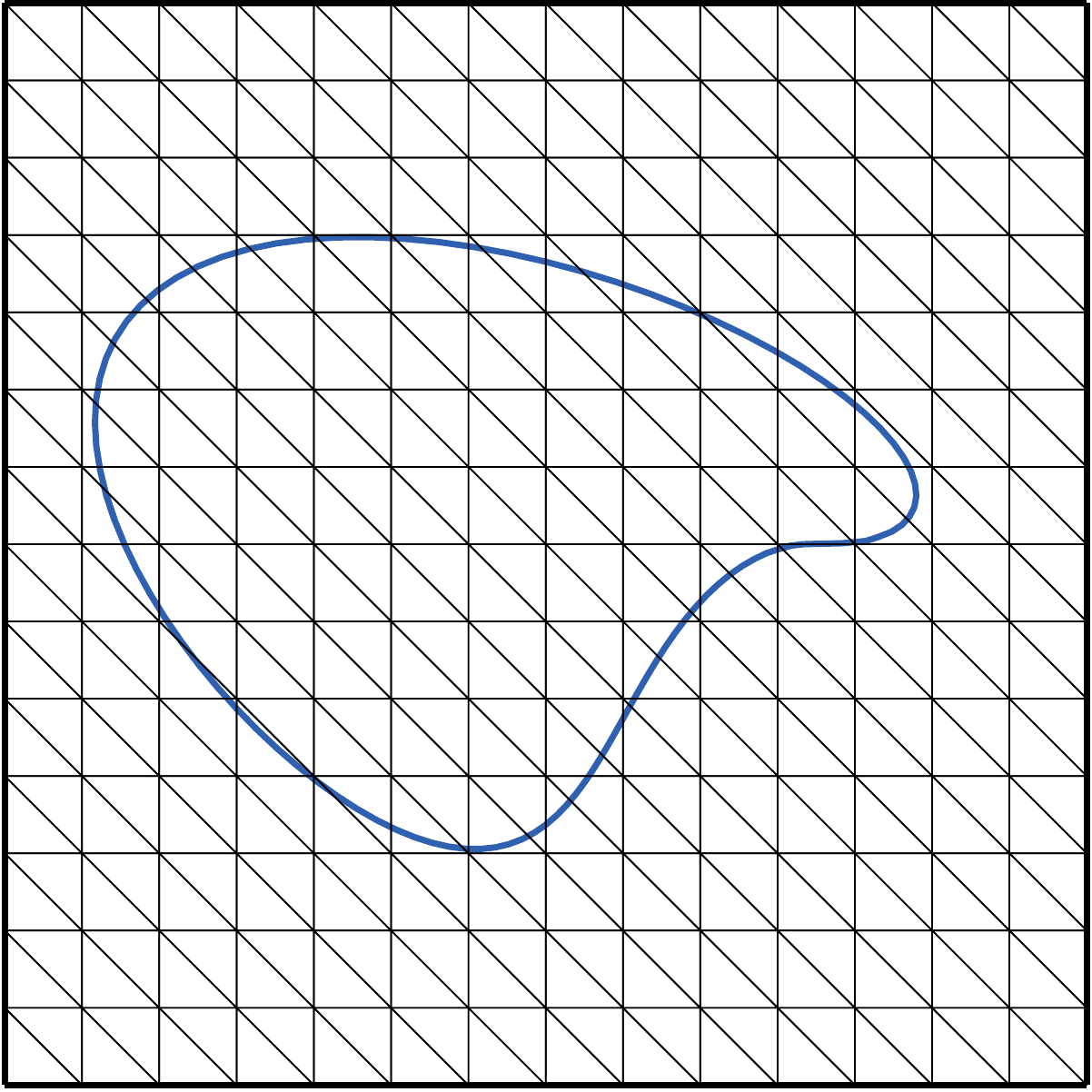} \hspace{1cm}
\includegraphics[scale=0.45]{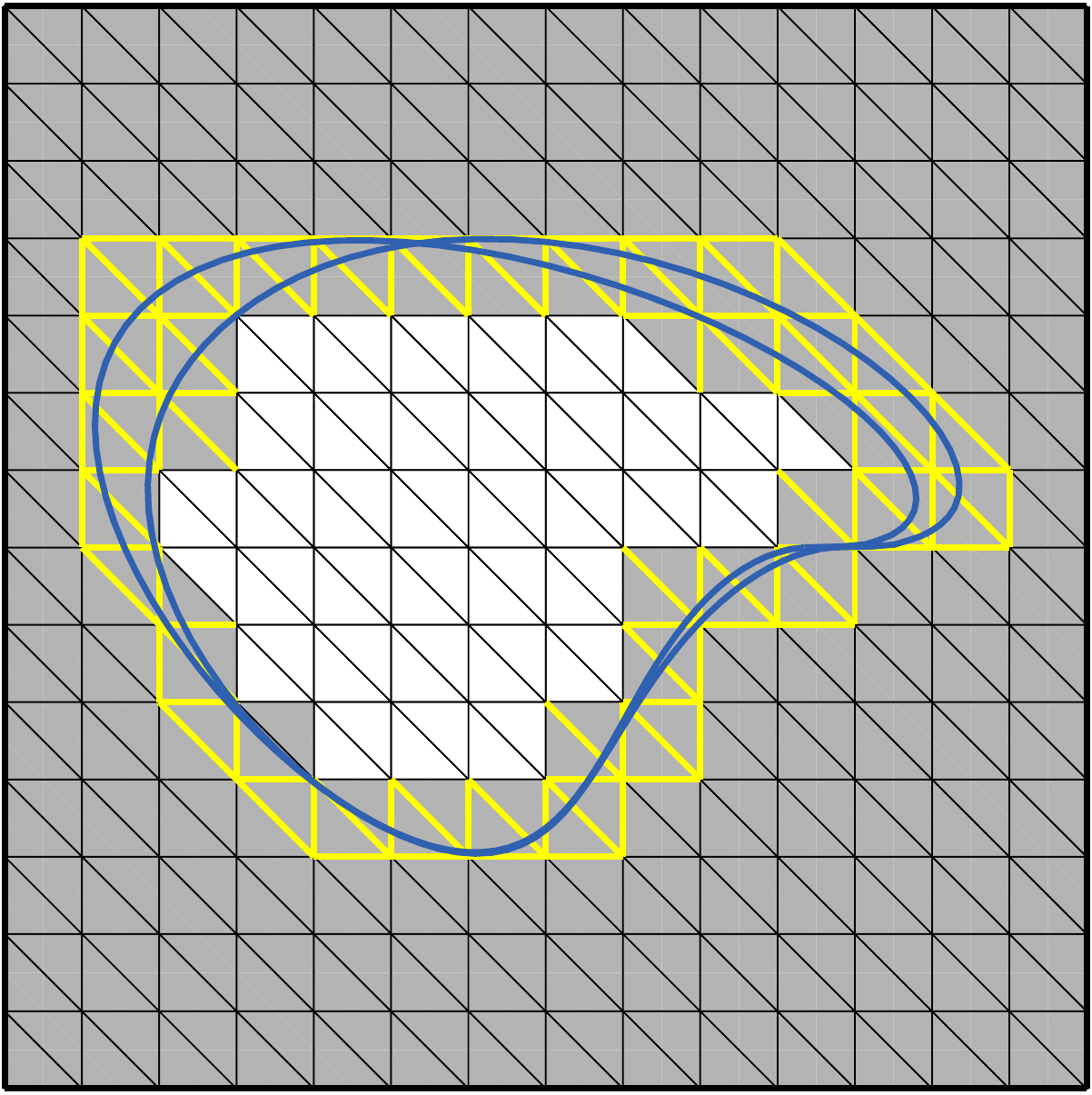}
\caption{Illustration of the sets introduced in Section~\ref{sec:mesh}. Left: a fixed background mesh of the computational domain $\Omega$. Right: the active mesh associated with the subdomain $\Omega_1$ and a time interval $I_n = (t_{n},t_{n+1}]$. The two blue curves show the position of the interface $\Gamma(t)$ at the endpoints $t=t_{n}$ and $t=t_{n+1}$ of the time interval $I_n$. The shaded domain shows $\mcN_{h,1}^n$ and the elements of the active mesh. Edges in $\mcF_{h,1}^n$ are marked with yellow thick lines.\label{fig:illustmesh}}.
\end{figure} 
Let $\mcK_h$ be a quasi-uniform partition of $\Omega$ into shape regular simplices of diameter $h$ and $\mcK_{h/2}$ the mesh obtained by refining $\mcK_h$ uniformly once. These meshes are time independent, we refer to them as fixed background meshes, they are generated independently of the location of the interface $\Gamma$ and thus the interface may at any time cut through these meshes arbitrarily. 

Let $0=t_0 < t_1 < \dots < t_N = T$ be a partition of $I$ into steps $I_n = (t_{n}, t_{n+1}]$ of length $\Delta t_n = t_{n+1} - t_{n} \ , \ n = 0,1,\dots, N-1$. For each time interval $I_n$ and subdomain $\Omega_i$ we define a so called active mesh which is the set of elements in the fixed background mesh that have a nonempty intersection with the subdomain $\Omega_i(t)$ for any time $t\in I_n$. More precisely, we define  
\begin{align}
	\mcN_{h,i}^n & = \bigcup_{t \in I_n} \bigcup_{K \in \mcK_{h,i}(t)} K \ , \ i=1,2, \\
	\mcK_{h,i}(t) &= \enstq{K \in \mcK_h}{ |\bar{\Omega}_i(t) \cap \partial K| > 0}, 
\end{align}
(where $| \cdot|>0$ means a positive surface measure) and the active mesh $\mcK_{h,i}^n$ is the set of elements that constitute the domain $\mcN_{h,i}^n$.
We also need the set of elements in the fixed background mesh that exhibit a nonempty intersection with the interface for any time $t\in I_n$:
\begin{align}
  \mcK_{h,\Gamma}^n =\enstq{K \in \mcK_{h} }{|\bar{K} \cap \Gamma(t) | > 0, t \in I_n}.
\end{align}
Finally, we denote the set of faces in $\mcK_{h,\Gamma}^n$ which are shared by two elements in the active mesh $\mcK_{h,i}^n$ by $\mathcal{F}_{h,i}^n$.
For an illustration, in two space dimensions,  of the sets we introduced in this section see Fig. \ref{fig:illustmesh}.

\subsection{Finite element spaces}\label{sec:space}
We take the inf-sup stable P1-iso-P2/P1 linear element pair as in \cite{HaLaZa14} but one may also choose other element pairs, see e.g. \cite{GuOl18}. Thus, on the fixed background mesh $\mcK_h$ we define 
\begin{equation}
 Q_{h}=\enstq{q_h \in C^0(\Omega)}{q_h\restreinta{K} \in P_1(K), \ \forall K \in \mcK_h, \ \prodscal{\mu^{-1}q_h}{1}_{\tdomain} = 0},
\end{equation}
the space of continuous piecewise linear functions with $\prodscal{\mu^{-1}q_h}{1}_{\tdomain} = 0$, and on the fixed background mesh $\mcK_{h/2}$ we let $V_{h/2}$ be the space of vector valued continuous piecewise linear functions.

On the space-time slab $I_n \times \mcN_{h,i}^n$ and  $I_n \times \mcN_{h/2,i}^n$ we define the spaces 
\begin{align}\label{eq:spaceQhi}
	Q_{h,i}^{n,r} & = P_r(I_n)  \otimes Q_{h}\restreinta{\mcN_{h,i}^n}  \quad i=1,2
\end{align}
and 
\begin{align}\label{eq:spaceVhi}
	V_{h/2,i}^{n,r} & = P_r(I_n)  \otimes V_{h/2}\restreinta{\mathcal{N}_{h/2,i}^n}  \quad i=1,2, 	
\end{align}
respectively. Here $P_{r}(I_n)$ is the space of polynomials of degree less or equal to $r$ on the interval $I_n$.  
We are now ready to define our pressure and velocity space:
\begin{align*}
 Q^{n,r}_h = & 
 	\enstq{q_h = (q_{h,1},q_{h,2})}{q_{h,i} \in Q^{n,r}_{h,i}, \ i=1,2}, \\ 
 V^{n,r}_{h} = & 
 \enstq{\bfv_h = (\bfv_{h,1},\bfv_{h,2})}{\bfv_{h,i} \in V^{n,r}_{h/2,i}, \ i=1,2}. 
\end{align*}

\begin{figure}
\centering
	\includegraphics{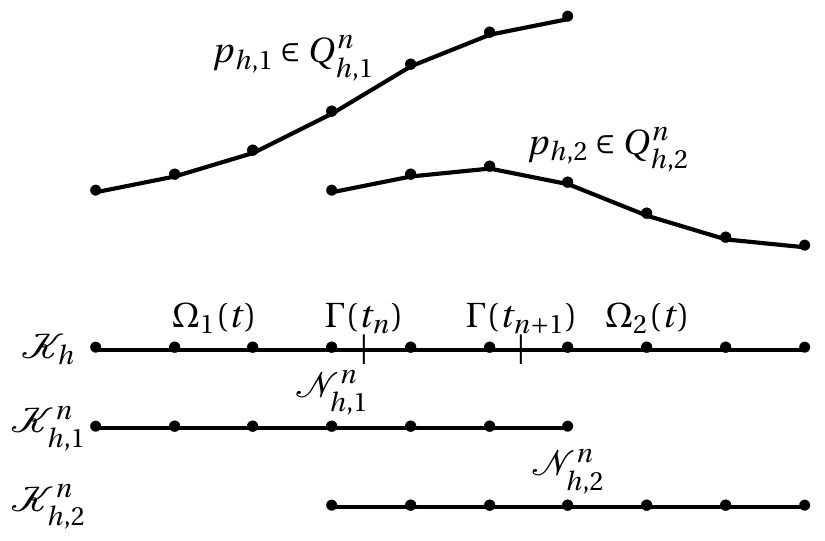}
	 \caption{Illustration of the active meshes and the pressure space associated with a time interval $I_n$ in a one space dimensional model. The interface evolves from $\Gamma(t_{n})$ to $\Gamma(t_{n+1})$ during the time interval $I_n$. \label{fig:illustcutfemspace}} 
\end{figure}

 Note that functions in the pressure space $Q^{n,r}_h$ and the velocity space $V^{n,r}_{h}$ consist of pair of functions associated to the two subdomains $\Omega_i$, $i=1,2$ and are double valued on elements in $\mcK_{h,\Gamma}^n$, since those elements exist in both active meshes $\mcK_{h,1}^n$ and $\mcK_{h,2}^n$. Hence functions in the pressure space as well as functions in the velocity space can be discontinuous at the interface $\Gamma(t)$ for $t\in I_n$. See Fig. \ref{fig:illustcutfemspace} for an illustration in one space dimension. 

A function  $\bfv_h \in V^{n,r}_h$ can be written in the following form
\begin{align*}
\bfv_h(t, \bfx) =  \left( \sum_{j=0}^r \bfv_{h,1,j} \left(\frac{t-t_{n}}{\Delta t_n} \right)^j, 
\sum_{j=0}^r \bfv_{h,2,j} \left(\frac{t-t_{n}}{\Delta t_n} \right)^j \right),
\end{align*}
where $t \in I_ n$ and for each j    
\[
	\bfv_{h,i,j} = \sum_{k=1}^{N_i} \textbf{$\xi$}_{i,k,j} \phi_k(\bfx)\restreinta{\mcN_{h,i}^n} \ , \ i=1,2. 
\]
Here $\textbf{$\xi$}_{i,k,j} \in \R^d$ are coefficients, $N_i$ is the number of nodes in the active mesh $\mcK_{h,i}^n$ and $\phi_k(\bfx)$ is the standard nodal basis function associated with mesh vertex $k$. 
To obtain a second order accurate discretization in time we propose to let the trial and the test functions be piecewise linear in time, i.e. $r=1$. Note that we then for $t\in I_n$ and $\bfu_h$ and $\bfv_h$ in $V^{n,1}_h$ have
\begin{equation}
 \prodscal{\rho \partial_t \bfu_h}{\bfv_h}_{\tdomain} = \frac{1}{\Delta t_n} \sum_{i=1}^2\prodscal{\rho_i\bfu_{h,i,1}}{\bfv_{h,i,0}}_{\Omega_i(t)} + \frac{t-t_{n}}{(\Delta t_n)^2}\sum_{i=1}^2 \prodscal{\rho_i\bfu_{h,i,1}}{\bfv_{h,i,1}}_{\Omega_i(t)}. 
\end{equation}

\subsection{The finite element method} \label{sec:weakform_nummeth}
We propose the following variational formulation where space and time are treated similarly and interface and boundary conditions are imposed weakly: for every time interval $I_n$ given $\bfu_h(t_{n}^-,\vect{x})$ find $(\bfu_h, p_h) \in V^{n,r}_h \times Q_h^n$ such that 
\begin{align}\label{eq:weakformcutfem}
A^n_h(\bfu_h,p_h,\bfv_h,q_h) +  S^n_{\bfu}( \bfu_h, \bfv_h) + S^n_p(p_h,q_h) = L^n_h(\bfv_h,q_h)
 \quad \forall (\bfv_h,q_h) \in V^{n,r}_h \times Q^{n,r}_h.
\end{align}
Here
\begin{align}\label{eq:weakformAh}
A^n_h(\bfu_h,p_h,\bfv_h,q_h)&=\int_{I_n} \prodscal{\rho \partial_t \bfu_h}{\bfv_h}_{\tdomain} \diff t  +
\int_{I_n} \prodscal{\rho(\bfu_h \cdot \grad){\bfu_h}}{\bfv_h}_{\tdomain} \diff t + 
\nonumber \\
+ & \int_{I_n} a(t, \bfu_h,\bfv_h)-b(t,\bfv_h,p_h)+b(t,\bfu_h,q_h) \diff t
+  \prodscal{\rho\jump{\bfu_h}_{t_n}}{\bfv_h(t^+_{n},\bfx)}_{\Omega_1(t_{n})\cup \Omega_2(t_{n})},
\end{align}
$\jump{\bfu_h}_t = \bfu_h(t^+, \bfx) - \bfu_h(t^-, \bfx)$, $t^{\pm} := \lim_{\varepsilon \to 0} t \pm \varepsilon$, 
\begin{align*}
	L^n_h(v_h,q_h) = \int_{I_n} l(t,\bfv_h,q_h) \diff t,
\end{align*}
the forms $a(t, \bfu,\bfv)$, $b(t,\bfv,p)$, and $l(t,\bfv,q)$ are defined in equation \eqref{eq:forma}-\eqref{eq:forml}. We state our choice for the weights $\w_1$, $\w_2$ and the penalty parameters $\lambda_\Gamma$ and $\lambda_{\partial \Omega}$ used in these forms in Section \ref{sec:parameters}. Note that via the last term in equation \eqref{eq:weakformAh} and the known solution $\bfu_h(t_{n}^-,\vect{x})$ from the previous space-time slab we weakly enforce continuity at time $t_{n}$. Starting from the initial condition $\bfu_h(t_{0}^-,\vect{x})=\bfu^0(\bfx)$ in $\Omega_1(t_{0})\cup \Omega_2(t_{0})$ we can solve the variational formulation \eqref{eq:weakformcutfem} one space-time slab at a time.   

The terms $S_p^n(p_h,q_h)$ and $S_{\bfu}^n(\bfu_h,\bfv_h)$ in equation \eqref{eq:weakformcutfem} are appropriate stabilization terms. These stabilization terms are used to 1) control the condition number of the resulting system matrix; 2) ensure the inf-sup stability of the method; 3) allow us to directly approximate the space-time integrals in the variational formulation using quadrature rules, see Section \ref{sec:quadintime}.
We use  
\begin{align*}
	S_p^n(p_h,q_h) = \int_{I_n} s_p(t,p_h,q_h) \diff t \ \text{ and } 
	\  S_{\bfu}^n(\bfu_h,\bfv_h) = \int_{I_n} s_{\bfu}(t,\bfu_h,\bfv_h) \diff t,
\end{align*}
with
\begin{align}\label{eq:sp}
s_p(t,p_h,q_h) = \sum_{i=1}^2 \sum_{F \in \mathcal{F}_{h,i}^n} 
C_p \mu_i^{-1} h^3 \prodscal{\jump{ \bfn_F \cdot \grad{p_{h,i}}}_F}{\jump{ \bfn_F \cdot \grad{q_{h,i}}}_F}_F
\end{align}
and
\begin{align}\label{eq:su}
s_{\bfu}(t,\bfu_h,\bfv_h) = \sum_{l=1}^d \sum_{i=1}^2 \sum_{F \in \mathcal{F}_{h,i}^n}
C_{\bfu} \mu_i h \prodscal{\jump{ \bfn_F \cdot \grad{\bfu_{h,i}^l}}_F}{\jump{ \bfn_F \cdot \grad{\bfv_{h,i}^l}}_F}_F.
\end{align}
Here $C_p$ and $C_{\bfu}$ are positive constants and $\jump{v}_F$ denotes the jump of a function $v$ at the face F and is defined as $\jump{v}_F=v^+-v^-$, where $v^\pm=\lim_{t\rightarrow 0^+} v(\bfx \mp t \bfn_F)$, $\bfx\in F$, and $\bfn_F$ is a fixed unit normal to $F$.
We usually choose the constants $C_p$ and $C_{\bfu}$ in the interval $\left[10^{-3},  10^{-1}\right]$. Increasing the constants $C_p$ and $C_{\bfu}$ typically leads to a smaller condition number  of the resulting linear system but larger error.
The specific form of the stabilization given in equation \eqref{eq:sp} and \eqref{eq:su} is for functions that are continuous and piecewise linear in space, for higher order elements in space one needs to include higher order derivatives or consider other type of stabilization terms, see Remark \ref{rem:highorder}. We also emphasize that for each space-time slab the set $\mathcal{F}_{h,i}^n$ is fixed and does not depend on $t$, however the trial and the test functions are time dependent.

\begin{remark}\label{rem:highorder}
Using the proposed discretization with $r=1$ (functions are piecewise linear in time) and provided that all space-time integrals can be approximated with at least second order accuracy we expect to get an approximate velocity field $\bfu_h(t,\bfx)$ which is second order accurate in both space and time. In Section \ref{sec:surften} we present how to obtain an accurate approximation of the surface tension force.
If the regularity of the problem allows, higher order accuracy than second order can be obtained with the presented space-time cut finite element strategy. In that case the P1-iso-P2/P1 linear element pair has to be changed to higher order elements, $r$ should increase, all integrals in both space and time have to be accurately approximated,  and the ghost penalty stabilization $s_{\bfu}$ and $s_p$ has to be updated. For example the following ghost penalty stabilization can be used 
\begin{equation}
s_{\bfu}(t,\bfu_h,\bfv_h) = \sum_{l=1}^d \sum_{i=1}^2 \sum_{F \in \mathcal{F}_{h,i}^n}\sum_{m=0}^M
C_{\bfu,m} (\mu_i,\rho_i) h^{2m-1} \prodscal{\jump{ D^m_{\bfn_F}{\bfu_{h,i}^l}}_F}{\jump{D^m_{\bfn_F} {\bfv_{h,i}^l}}_F}_F,
\end{equation}
where $\jump{D^m_{\bfn_F} v}_F$ denotes the jump in the normal derivative of order m across the face F and $C_{\bfu,m} (\mu_i,\rho_i)$ is a constant that may depend on $\mu_i$, $\rho_i$, and the order m of the directional derivative. Here $\bfu_{h,i}$ and $\bfv_{h,i}$ are piecewise polynomials of degree M in space. The stabilization $s_p$ should be updated similarly to include higher order derivatives.  
\end{remark}

\subsubsection{Quadrature in time} \label{sec:quadintime}
In the proposed space-time method as in \cite{HanLarZah16, Zah18} we don't explicitly construct the space-time domain in $\R^{d+1}$. Due to the added stabilization all space-time integrals in the variational formulation \eqref{eq:weakformcutfem} can be approximated using quadrature rules, first in time and then in space. Thus, given an accurate quadrature rule in the time interval $I_n$ with $N_m$ weights ($\omega_m^n$) and quadrature points ($t_m^n$) the discrete formulation is:  given the solution $\bfu_h(t^-_{n},\bfx)$ from the previous space-time slab find
$(\bfu_h, p_h)\in  V^{n,r}_h \times Q^{n,r}_h$ such that 
\begin{align}\label{eq:weakformquad}
&\sum_{m=1}^{N_m} \omega_m^n \left(
\prodscal{\rho \partial_t \bfu_h}{\bfv_h}_{\Omega_1(t_m^n) \cup \Omega_2(t_m^n)} +
 \prodscal{\rho(\bfu_h \cdot \grad){\bfu_h}}{\bfv_h}_{\Omega_1(t_m^n) \cup \Omega_2(t_m^n)} \right) + 
\nonumber \\
+ & \sum_{m=1}^{N_m} \omega_m^n  a(t_m^n, \bfu_h,\bfv_h)-b(t_m^n,\bfv_h,p_h)+b(t_m^n,\bfu_h,q_h) 
+  \prodscal{\rho \bfu_h(t^+_{n},\bfx)}{\bfv_h(t^+_{n},\bfx)}_{\Omega_1(t_{n})\cup \Omega_2(t_{n})}
\nonumber \\
+ &
\sum_{m=1}^{N_m} \omega_m^n \left(s_p(t_m^n,p_h,q_h) + s_{\bfu}(t_m^n,\bfu_h,\bfv_h) \right) =
\nonumber \\
=& \sum_{m=1}^{N_m} \omega_m^n l(t_m^n,\bfv_h,q_h) +
\prodscal{\rho \bfu_h(t^-_{n},\bfx)}{\bfv_h(t^+_{n},\bfx)}_{\Omega_1(t_{n})\cup \Omega_2(t_{n})}
\end{align}
for all $(\bfv_h, q_h) \in V^{n,r}_h \times Q^{n,r}_h$.

To have a second order accurate discretization in time we need $r=1$ and a quadrature rule which is at least second order, i.e. has a degree of precision of at least 1.
Numerical examples in \cite{HanLarZah16, Zah18} using $r=1$ with both the trapezoidal rule and the Simpson's rule show second order convergence. Note that these quadrature rules,  see Table \eqref{tab:quadrule}, and in general closed Newton-Cotes formulas include the endpoints of the time interval $I_n$ and therefore some computations can be reused when passing from one space-time slab to another. However, we emphasize that other accurate quadrature rules can also be used. 
\begin{table}[h]
\centering
\begin{tabular}{|l|l|l|l|}
\toprule
 $N_m$ & quadrature points  $t_m^n$ & quadrature  weights $\omega_m^n$  & degree of precision \\  \midrule
 $2$ & $t_1^n=t_{n-1}$, $t_2^n=t_{n}$ &$\omega_1^n=\omega_2^n=\frac{k_n}{2}$  & 1\\
 &  & &\\
$3$ & $t_1^n=t_{n-1}$,  $t_3^n=t_{n}$, $t_2^n=\frac{t_{n-1}+t_n}{2}$ & $\omega_1^n=\omega_3^n=\frac{k_n}{6}$, $\omega_2^n=\frac{4k_n}{6}$ & 3\\
\bottomrule
\end{tabular}
\caption{First row: Trapezoidal rule. Second row: Simpson's rule. \label{tab:quadrule}}
\end{table}
 
Note that the stabilization terms are integrated on faces in $\mathcal{F}_{h,i}^n$, independent of time $t\in I_n$. Since the trial and the test functions are both polynomials of degree $r$ in time these stabilization terms are polynomials of degree $2r$ in time and the integration in time over $I_n$ can be done analytically. We find the faces in the set $\mathcal{F}_{h,i}^n$ in the following way:  1) We find the position of the interface $\Gamma(t)$ and the domains $\Omega_i(t)$ at the discrete time instances $T_n=\{t_n, \{t_m^n\}_{m=1}^{N_m}, t_{n+1}\}$. 2) We define an element $K$ to be in $\mcK_{h,\Gamma}^n$ if it is cut by the interface for some time $t \in T_n$ or if there are two time instances $t_k \in T_n$ and $t_l \in T_n$ such that $K$ is in $\Omega_1(t_k)$ but not in $\Omega_1(t_l)$ (i.e. $K\in \Omega_2(t_l)$). If for example a signed distance function is available at the time instances $t \in T_n$ this information can be determined from the sign of that function. All faces in $\mcK_{h,\Gamma}^n$ except those that are in $\Omega_2(t)$ for all $t\in T_n$ are in $\mathcal{F}_{h,1}^n$ and are stabilized. Faces in $\mathcal{F}_{h,2}^n$ are defined similarly. 
In step 1 note that, from the previous space-time slab, we already know the position of the interface and the subdomains at time $t=t_n$. Therefore for each $I_n$ we find the interface and the subdomains at $N_{m}-1$ time instances, if the quadrature rule does include the endpoints $t_n$ and $t_{n+1}$, otherwise at $N_{m}+1$ points. For example using the trapezoidal rule we have to find the interface  at one time instance $t=t_{n+1}$ in each interval $I_n$ while if we use the Simpson's rule we find the interface at two time instances $t_{n+1/2}$ and $t_{n+1}$  in each interval $I_n$.

\begin{remark}
Letting the trial and the test functions be piecewise constant in time, i.e. choosing $r=0$, one can obtain a scheme that corresponds  to using backward Euler for the time discretization and CutFEM in space. Take $r=0$ in equation \eqref{eq:spaceQhi} and \eqref{eq:spaceVhi}. The term $\partial_t \bfu_h$ vanishes. Let  $(\bfu_h^{n-1}, p_h^{n-1})$ denote the solution $(\bfu_h(t,\bfx), p_h(t,\bfx))$ from the previous space-time slab, which is constant in time for $t\in I_{n-1}$ and let  $(\bfu_h^{n}, p_h^{n})$ be the solution on the current space-time slab, i.e. $t\in I_{n}$, we have that the last term in equation \eqref{eq:weakformAh} is 
\begin{equation}
\prodscal{\rho\jump{\bfu_h}_{t_n}}{\bfv_h(t^+_{n},\bfx)}_{\Omega_1(t_{n})\cup \Omega_2(t_{n})}=\prodscal{\rho \left(\bfu_h^{n}(\bfx)-\bfu_h^{n-1}(\bfx)\right)}{\bfv_h(\bfx)}_{\Omega_1(t_{n})\cup \Omega_2(t_{n})}.
\end{equation}
If we in equation \eqref{eq:weakformquad} now use a quadrature rule in the time interval $I_n$ defined by $N_m=1$, weight $\omega_1^n=\Delta t_n$, and quadrature point $t_1^n=t_{n}$ we have the discrete formulation: given $\bfu_h^{n-1}$ find $(\bfu_h^{n}, p_h^{n})\in  V^{n,0}_h \times Q^{n,0}_h$ such that 
\begin{align}\label{eq:weakformEuler}
& \prodscal{\rho \bfu_h^{n}(\bfx)}{\bfv_h(\bfx)}_{\Omega_1(t_{n})\cup \Omega_2(t_{n})}
+
\Delta t_n\prodscal{\rho(\bfu_h^{n} \cdot \grad){\bfu_h^{n}}}{\bfv_h}_{\Omega_1(t_{n}) \cup \Omega_2(t_n)}  + 
\nonumber \\
+ &
\Delta t_n  a(t_{n}, \bfu_h^{n},\bfv_h)-\Delta t_n b(t_{n},\bfv_h,p_h^{n})+\Delta t_n b(t_{n},\bfu_h^{n},q_h)
+  
\nonumber \\
+ &
 \Delta t_n  s_p(t_{n},p_h^{n},q_h) +  \Delta t_n s_{\bfu}(t_{n},\bfu_h^{n},\bfv_h) =
\nonumber \\
=&  \Delta t_n l(t_{n},\bfv_h,q_h) +
\prodscal{\rho \bfu_h^{n-1}(\bfx)}{\bfv_h(\bfx)}_{\Omega_1(t_{n})\cup \Omega_2(t_{n})}
\end{align}
for all $(\bfv_h, q_h) \in V^{n,0}_h \times Q^{n,0}_h$. This is the same discrete formulation one would obtain using backward Euler for the time discretization and a CutFEM in space. 
\end{remark}

\begin{remark}
All computations, including the construction of the interface $\Gamma(t)$ and the subdomains $\Omega_1(t)$ and $\Omega_2(t)$ are only done at discrete time instances. We never explicitly construct the space-time domain in $\R^{d+1}$. Given a method for representing and evolving the interface it is straightforward to implement the proposed space-time CutFEM from a stationary CutFEM.  For space-time methods which are built on explicitly constructing the space-time domain see e.g. \cite{CL15}. 
\end{remark}

\subsubsection{The weights and the penalty parameters}\label{sec:parameters}
There are several choices for the weights  $\w_1$ and $\w_2$ in the averaging operators \eqref{eq:averagop} and for the penalty parameters $\lambda_\Gamma$ and $ \lambda_{\partial \Omega }$ in the Nitsche penalty terms, see e.g. \cite{BuZu11, AnnHauDol12, WZKB, HaLaZa14}. However, the weights should be chosen so that $\w_1+\w_2=1$. The penalty parameters are as $C_{p}/h$, where $C_{p}$ is a sufficiently large constant. In \cite{HaLaZa14}, based on the analysis, we suggested to choose these parameters locally as 
\begin{equation}
%\[
	\w_1 \restreinta{K}= \frac{\mu_2 \alpha_{1,K}}{\mu_1\alpha_{2,K} + \mu_2 \alpha_{1,K}} \ , \
	\w_2  \restreinta{K}= \frac{\mu_1 \alpha_{2,K}}{\mu_1\alpha_{2,K} + \mu_2 \alpha_{1,K}},
%\]
\end{equation}
and 
\begin{equation}
%\[
	\lambda_{\Gamma(t)} \restreinta{K} = \frac{\{\mu \}}{h_K} 
	\left( D + C \frac{\gamma_K}{\alpha_K} \right)
	\ , \
%\]
\end{equation}
where $K$ is an element cut by the interface with $| K \cap \Omega_i(t)| = \alpha_{i,K} h_K^d$, $\alpha_K = \alpha_{1,K}+\alpha_{2,K}$, and $|K \cap \Gamma(t)| = \gamma_K h_K^{d-1}$. Here $C>1$ and $D>0$ are constants. For an element $K$ on the boundary we defined 
$$\lambda_{\partial \Omega } \restreinta{K \cap \Omega_i} = 
		\frac{\mu_i}{h_K} \left(G + H \frac{\gamma_{\partial \Omega , K}}{\alpha_K} \right)$$
with $|K \cap \partial \Omega| = \gamma_{\partial \Omega , K} h_K^{d-1}$, and $G>0$, $H>0$ sufficiently large \cite{HaLaZa14}. Under the assumption that the interface does not cut the boundary of the domain $\Omega$  the stabilization term $s_{\bfu}$ could with this choice of parameters be chosen weaker (an $h^3$ scaling in \eqref{eq:su}) and coercivity of $a(\cdot , \cdot)$ would still be ensured, see Remark 1 in \cite{HaLaZa14}. With this choice of parameters the size of the penalty parameters are minimized. However, if the stabilization term $s_{\bfu}$ is chosen as in \eqref{eq:su} (with an $h^1$ scaling) the weights do not need to include the scaling with the relative area/volume of each subdomain,  $\alpha_{i,K}$, and can be chosen as 
\begin{equation}
\w_1=\frac{\mu_2 }{\mu_1 + \mu_2} \ , \
\w_2 = \frac{\mu_1}{\mu_1+ \mu_2}.
\end{equation}
Another choice taking into account high contrasts in both viscosity and the density (for $\rho_i>0$) is thus 
\begin{equation}
\w_1=\frac{\mu_2/\rho_2 }{\mu_1/\rho_1 + \mu_2/\rho_2} \ , \
\w_2 = \frac{\mu_1/\rho_1}{\mu_1/\rho_1 + \mu_2/\rho_2}.
\end{equation}
Note that in all cases $\w_1+\w_2=1$.

\section{The representation and evolution of the interface} \label{section - level set}
In this work we use the level set method \cite{OshFed01, Set01} to represent and evolve the interface. However, we emphasize that the finite element method presented in the previous section is independent of the numerical technique used for representing and evolving the interface and other methods such as e.g. a front-tracking method \cite{TryBuEsetal01} can also be used. 

Let $\lsf(t,\bfx) : I \times \R^d \rightarrow \R$ be the signed distance function with positive sign in $\Omega_2(t)$, the subdomain enclosed by the interface. The zero level set of this function represents the interface $\Gamma(t)$. 
The unit normal is defined as $\bfn(t,\bfx) = \grad{\lsf}(t,\bfx)$ for $\bfx \in \Gamma(t)$. 
Given a vector field $\vel \in \R^d$ and $t\in I$, the evolution of the interface $\Gamma(t)$ is governed by the following partial differential equation 
\begin{align}
	\partial_t \lsf + \vel \cdot \grad{\lsf} = 0  \quad \textrm{in } \Omega \label{eq : Level set SF}
\end{align}
with initial condition $\lsf(0, \bfx) = \lsf_0(\bfx) $ given by the initial configuration of $\Gamma$. 

Denote by $W_{h/2,q}$ the space of continuous piecewise polynomials of degree less than or equal to $q \geq 1$ defined on the fixed background mesh $\mcK_{h/2}$. We denote by $\lsf_{h,q}$ an approximation of the level set function in $W_{h/2,q}$. The continuous piecewise linear approximation of $\lsf$ on $\mcK_{h/2}$ is 
\begin{equation}\label{eq:p1repoflsf}
\lsf_{h,1}=I_h^1 \lsf_{h,q} \in W_{h/2,1},
\end{equation}
 where $I_h^1$ is the nodal interpolation operator on $W_{h/2,1}$. 

We discretize \eqref{eq : Level set SF} using the Crank-Nicolson scheme and quadratic elements in space with a streamline diffusion stabilization: given $\lsf_{h,2}^{k-1} \in W_{h/2,2}$ find $\lsf_{h,2}^k \in W_{h/2,2}$, such that 
\begin{align}\label{eq:advnum}
\prodscal{\frac{\lsf_{h,2}^k}{\Delta t_k}  + \frac{1}{2} \vel^k \cdot \nabla \lsf_{h,2}^k}{v}_{\tdomain} 
+ \prodscal{\frac{\lsf_{h,2}^k}{\Delta t_k} + \frac{1}{2} \vel^k \cdot \nabla \lsf_{h,2}^k}{\tau_{SD}\vel^k \cdot \nabla v}_{\tdomain}  = \nonumber \\
\prodscal{\frac{\lsf_{h,2}^{k-1}}{\Delta t_k}  - \frac{1}{2} \vel^{k-1} \cdot \nabla \lsf_{h,2}^{k-1}}{v}_{\tdomain} 
+ \prodscal{\frac{\lsf_{h,2}^{k-1}}{\Delta t_k} \lsf_{h,2}^{k-1} - \frac{1}{2} \vel^{k-1} \cdot \nabla \lsf_{h,2}^{k-1}}{\tau_{SD}\vel^k \cdot \nabla v}_{\tdomain}  
\end{align}
for all $v \in W_{h/2,2}$ and where $\tau_{SD} = 2\left( \Delta t_k^{-2} + |\vel|^2 h^{-2} \right)^{-\frac{1}{2}}$ is the streamline diffusion parameter. Note that 
$\lsf_{h,2}^k$ is a piecewise quadratic approximation of the signed distance function on the refined mesh $\mcK_{h/2}$ at time instance $t_k$.
The time instances we use here are associated to the quadrature points used in time to approximate the space-time integrals in the proposed finite element method, see Section \ref{sec:quadintime}. For example if the trapezoidal rule is used the points $t_k$ are exactly the endpoints of the time intervals $I_n$ and $\Delta t_k=\Delta t_n$. However, if Simpson's rule is used $\Delta t_k=1/2 \Delta t_n$ and $\{t_k\}$ include the endpoints and the midpoints of each interval $I_n$.

We take $\vel$ to be the fluid velocity, i.e. $\vel=\bfu_h$. 
Note that the integrals in \eqref{eq:advnum} are split into integrals over the subdomains and that $\bfu_h=\bfu_{h,i}$ in subdomain $\Omega_i$, see the previous section. 

In order to maintain a signed distance function  we discretize the following reinitialization equation 
\begin{align}\label{eq : Reinitialization SF}
&\partial_{\hat{t}} \rf = \sign{\lsf_h^k}(1 - |\grad{\rf}| ), \\ 
&\rf \restreinta{\hat{t}=0} = \lsf_h^k, \nonumber
\end{align}
as proposed in \cite{SuFa99}. The steady state solution of this problem yields a signed distance function but in practice, only a few steps are needed in order to obtain a signed distance function in a neighborhood around the interface.

To evaluate the integrals in the variational formulation presented in Section \ref{sec:nummeth}  we find the interface explicitly as the zero level set of the piecewise linear approximation $\lsf_{h,1}^k$ of the signed distance function, i.e. We find $\Gamma_{h,1}(t_k)=\enstq{\bfx \in \Omega}{\lsf_{h,1}^k(\bfx)=0}$ with $\lsf_{h,1}$ as in \eqref{eq:p1repoflsf} and $q=2$. 
Note that $\Gamma_{h,1}(t_k)$ is planar on each element in $\mcK_{h/2}$ which is cut by $\Gamma_{h,1}(t_k)$ and in three space dimensions it consists of triangles and quadrilaterals which can be subdivided into triangles. Thus, almost all the integrals in the proposed variational formulation can easily be computed with second order accuracy, see also \cite{MiGi07}. The surface tension force, is the only term in our weak formulation that we need a special treatment for in order to compute it with second order accuracy. We treat the approximation of this term in the next section.

\section{Surface tension force} \label{sec:surften}
Recall the weak formulation from Section \ref{sec:weakform_nummeth}, due to the normal stress jump condition, we have the following term
\begin{equation}\label{eq:surftenforce}
\prodscal{\sigma \kappa \bfn}{\langle \bfv_h \rangle}_{\Gamma(t)}, \quad \bfv_h \in V_h^{n,r}
\end{equation} 
in the form $l$, see equation \eqref{eq:forml}. We now present a numerical technique for computing a second order accurate (in the $L^2$-norm) mean curvature vector $\meancurv= \kappa \bfn$ to be used in the term \eqref{eq:surftenforce}.  

First note that for $t\in I$, 
 the mean curvature vector $\meancurv: \Gamma(t) \rightarrow \R^d$
satisfies the following weak problem:
find $\meancurv \in [H^1(\Gamma(t))]^d$ such that 
\begin{align}
	\prodscal{\meancurv}{\bfw}_{\Gamma(t)} 
	= \prodscal{\nabla_{\Gamma} \bfx_{\Gamma}}{\nabla_{\Gamma} \bfw}_{\Gamma(t)}, \quad 
\forall \bfw \in [H^1(\Gamma(t))]^d.
\label{eq : Mean curvature WF}
\end{align}
Here, $\bfx_\Gamma : \Gamma \ni \bfx \mapsto \bfx \in \R^d$ is the coordinate map, $\nabla_\Gamma$ is the tangential gradient which we define as $\nabla_\Gamma = P_\Gamma \nabla$, with $P_\Gamma = \bfI - \bfn \otimes \bfn$, and 
$\nabla_{\Gamma}\bfw = \bfw \otimes \nabla_\Gamma$ for a vector valued function $\bfw$. Note that to arrive at equation \eqref{eq : Mean curvature WF} the definition of the mean curvature vector in terms of the Laplace-Beltrami operator was used together with integration by parts \cite{Dz88}. In several works see e.g. \cite{Ba01, Hys06, GrRe07} the following form 
\begin{equation}\label{eq:surftenforceLB}
\prodscal{\sigma \kappa \bfn}{\bfv_h}_{\Gamma(t)}=
\prodscal{\sigma\nabla_{\Gamma} \bfx_{\Gamma}}{\nabla_{\Gamma} \bfv_h}_{\Gamma(t)}, 
 \quad \bfv_h \in V_h^{n,r}
\end{equation} 
is used in the discretization of the surface tension force and in this way the order of differentiation associated with the curvature is reduced. 

In this work, we propose to first compute a mean curvature vector $\meancurv_h$ based on the discrete Laplace-Beltrami operator and stabilization of the $L^2$ projection involved. Then use this stabilized mean curvature vector in the term \eqref{eq:surftenforce} as an approximation to $\kappa \bfn$. 
In \cite{HanLarZah15} we proposed a finite element formulation for computing such a stabilized  mean curvature vector from piecewise linear approximations of smooth surfaces and proved that $\meancurv_h$ would be a first order accurate approximation of the mean curvature vector in $L^2$. Note that in general no order of convergence in $L^2$ can be expected when computing a discrete mean curvature vector from a piecewise linear approximation of the interface. We now use the idea in \cite{HanLarZah15} to stabilize the $L^2$ projection and present a method for finding a $q$th order accurate approximation in $L^2$ of the mean curvature vector and consequently the surface tension force given an approximation of the interface as a  piecewise polynomial surface of order $q$.

\begin{figure}
\centering
\includegraphics[scale=0.45]{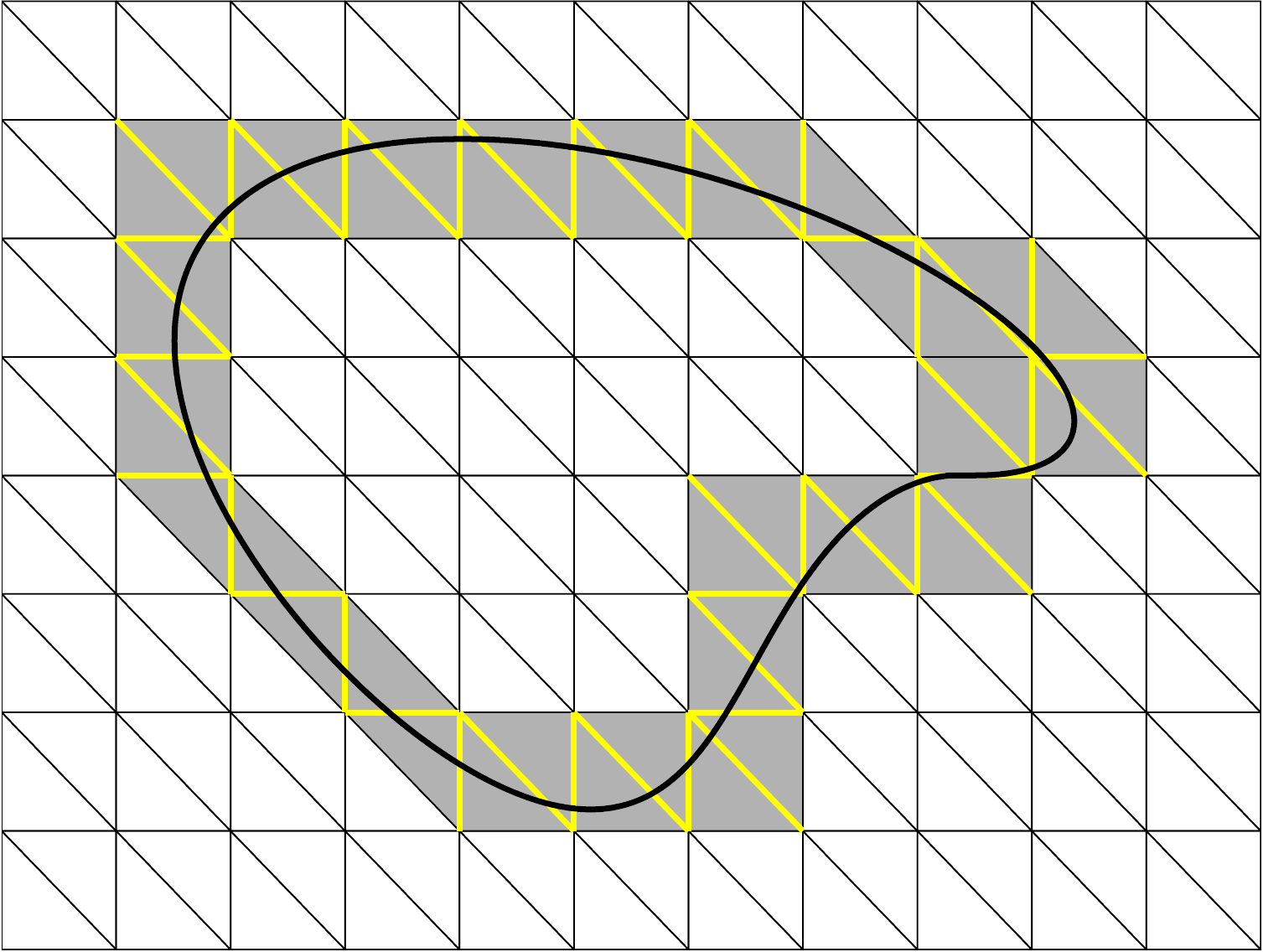}
\caption{Illustration of the active mesh $\mcK_{h, q,\Gamma}$ and the set $\mcF_{h,q,\Gamma}$. At a time $t\in I$, the shaded domain $\mcN_{h,q,\Gamma}(t)$ consist of elements in $\mcK_{h,q,\Gamma}(t)$ and the marked yellow edges show edges in $\mcF_{h,q,\Gamma}(t)$.  
\label{fig:illustactivemeshsurf}}.
\end{figure} 

For $t\in I$, let $\Gamma_{h,q}(t)$ be a piecewise polynomial surface of order $q \geq 1$ which is a $(q+1)$th-order accurate approximation of the interface $\Gamma(t)$. 
Let $W_{h,m}$ be the space of continuous piecewise polynomials of degree less than or equal to $m > 0$ defined on the fixed background mesh $\mcK_{h/2}$. 
On $\mcK_{h/2}$ we define the active mesh and the corresponding domain
\begin{align}
\mcK_{h,q,\Gamma}(t)&=\enstq{K \in \mcK_{h/2}}{|\bar{K} \cap \Gamma_{h,q}(t) | > 0}, \nonumber  \\
\mcN_{h,q,\Gamma}(t)&=\bigcup_{K\in\mcK_{h,q,\Gamma}(t)} K.
\end{align}
We now define the space $W_{h,m,q}(t)=W_{h,m}\restreinta{\mcN_{h,q,\Gamma}(t)}$.
Denote by $\mathcal{F}_{h,q,\Gamma}(t)$ the set consisting of internal faces (faces with two neighbors) in the active mesh, i.e. the set of all faces that are cut by the surface $\Gamma_{h,q}(t)$. See Fig. \ref{fig:illustactivemeshsurf} for an illustration of these sets in two space dimensions at some time $t$.
We define the following problem for the stabilized discrete mean curvature vector $\meancurv_h$:  find $\meancurv_h \in [W_{h,m,q}]^d$ such that
\begin{alignat}{2}\label{eq:numsurften}	
	\prodscal{\meancurv_h}{\bfw_h}_{\Gamma_{h,q}} + S_h(\meancurv_h, \bfw_h)
	= \prodscal{\nabla_{\Gamma_{h,q}} \bfx_{\Gamma_{h,q}}}{\nabla_{\Gamma_{h,q}} \bfw_h}_{\Gamma_{h,q}},	
	& & & \quad \forall \bfw_h \in [W_{h,m,q}]^d
\end{alignat}	
Here $\bfx_{\Gamma_{h,q}}: \Gamma_{h,q} \ni \bfx \mapsto \bfx \in \R^d$ is the discrete coordinate map, $\nabla_{\Gamma_{h,q}} = P_{\Gamma_{h,q}} \nabla$, is the tangential gradient with $P_{\Gamma_{h,q}} = \bfI - \bfn_{h,q} \otimes \bfn_{h,q}$ and $\bfn_{h,q}$ a $q$th-order accurate approximation of the interface normal  $\bfn$,
 and $S_h$ is an appropriate stabilization term. Note that $\nabla_{\Gamma_{h,q}} \bfx_{\Gamma_{h,q}}= P_{\Gamma_{h,q}}$.  We propose to choose 
\begin{equation}
m=q.
\end{equation}
Without the stabilization $S_h$ in equation \eqref{eq:numsurften} one expects a convergence of order $q-1$ in $L^2$ but with an appropriate stabilization term numerical experiments show convergence of order $q$ in  $L^2$, see Fig. \ref{fig:curvatureError}.
Here we use the stabilization proposed in \cite{Zah18, LarZah17} 
\begin{align}
	S_h(\meancurv_h, \bfw_h) = 
		 \sum_{j=1}^m \left(c_{F,j}
		 h^{2j-2} 
		 \prodscal{\jump{ D^j_{\bfn_F} \meancurv_h}_F}{\jump{ D^j_{\bfn_F} \bfw_h}_F}_{\mcF_{h,q,\Gamma}} 
		 +
		c_{\Gamma,j} h^{2j-2} 
		\prodscal{D^j_{\bfn_{{h,q}}} \meancurv_h}{ D^j_{\bfn_{{h,q}}} \bfw_h}_{\Gamma_{h,q}} \right),
\end{align}
where $c_{F,j}$ and $c_{\Gamma,j}$ are positive constants and $D^j_{\bfn}$ denotes the $j$th order directional derivative in the direction of $\bfn$.  This stabilization provides  control of the condition number both when linear as well as higher order elements are used \cite{LarZah17}  and improves the accuracy of the computed mean curvature vector. We emphasize that everything could also have been defined on the background mesh $\mcK_h$.

If the interface $\Gamma_{h,q}(t)$ is explicitly given standard quadrature rules can be used to compute the integrals in \eqref{eq:numsurften} and thus it is obvious how to compute all the integrals with high accuracy. However, in a level set method the piecewise polynomial surface $\Gamma_{h,q}(t)$ is implicitly defined as the zero level set of the level set function $\lsf_{h,q} \in W_{h/2,q}$. Several strategies for obtaining high order approximations of integrals on implicitly defined domains exist, see e.g. \cite{MuKuOb13, Sa15, FriOmer16, Le16}. 
We use the strategy in \cite{Le16} to accurately compute the integrals in \eqref{eq:numsurften}. This strategy which we describe in the next section, avoids integration on $\Gamma_{h,q}$ by a transformation of integrals on $\Gamma_{h,q}$ to integrals on the piecewise linear approximation $\Gamma_{h,1}$.

\subsection{Integration on implicitly defined interfaces} \label{section - isoparametric mapping}
We have that $\Gamma_{h,q}(t)$ is implicitly defined as the zero level set of the level set function $\lsf_{h,q} \in W_{h/2,q}$ and $q\geq 1$.
Following \cite{Le16} we introduce a transformation $\theta_h \in  [W_{h,q,\Gamma}]^d$ of the underlying mesh,  which maps the piecewise linear representation of the interface onto the zero level set of a high order approximation of the level set function, i.e. $\Gamma_{h,q}=\theta_h(\Gamma_{h,1})$.
Using this transformation we define the space $\hat{W}_{h,m,1} = W_{h,m,1} \circ \theta_h^{-1}$ 
and transform the weak formulation in equation \eqref{eq:numsurften} to: find $\meancurv_h \in [\hat{W}_{h,m,1}]^d$ such that
\begin{alignat}{2}\label{eq:numsurftenmapp}	
	\prodscal{\meancurv_h}{\bfw_h}_{\Gamma_{h,q}} + S_h(\meancurv_h, \bfw_h)
	= \prodscal{\nabla_{\Gamma_{h,q}} \bfx_{\Gamma_{h,q}}}{\nabla_{\Gamma_{h,q}} \bfw_h}_{\Gamma_{h,q}}	
	& & & \quad \forall \bfw_h \in [\hat{W}_{h,m,1}]^d.
\end{alignat}	

With this mapping, numerical integration on the implicitly defined surface $\Gamma_{h,q}$ can be implemented in the following way
\begin{align*}
\prodscal{\meancurv_h}{\bfw_h}_{\Gamma_{h,q}} 
& =  
\int_{\Gamma_{h,1}} \meancurv_h(\mapp_h(\bfx)) \cdot \bfw(\mapp_h(\bfx))
 	| \det(D\mapp_h(\bfx)) | \diff s_{h}  =
\int_{\Gamma_{h,1}} \tilde{\meancurv_h} \cdot \tilde{\bfw}
 	| \det(D\mapp_h(\bfx)) | \diff s_h 
\end{align*} 
and 
\begin{align*}
\prodscal{\nabla_{\Gamma_{h,q}} \bfx_{\Gamma_{h,q}}}{\nabla_{\Gamma_{h,q}} \bfw_h}_{\Gamma_{h,q}}
& =  
\int_{\Gamma_{h,1}} 
P_{\Gamma_{h,q}}(\mapp_h(\bfx))
:
P_{\Gamma_{h,q}}(\mapp_h(\bfx)) D\mapp_h^{-T}\grad{\tilde{\bfw}}
|\det(D\mapp_h(\bfx))| \cdot \norm{\mathbf{N}} \diff s_h.  	
\end{align*}
Here, $\tilde{\meancurv_h}=\meancurv_h(\mapp_h(\bfx)) \in [W_{h,m,1}]^d$, $\tilde{\bfw}=\bfw(\mapp_h(\bfx)) \in [W_{h,m,1}]^d$, $P_{\Gamma_{h,q}}(\mapp_h(\bfx)) = \bfI - \tilde{\bfn} \otimes \tilde{\bfn}$ with $\tilde{\bfn} = \mathbf{N} / \norm{\mathbf{N}}$,  $\mathbf{N} = (D\mapp_h)^{-T} \bfn_{h,1}$, and $\bfn_{h,1} = \grad{\phi_{h,1}} / \norm{\grad{\phi_{h,1}}}$.
The stabilization term is transformed in the same way. Thus, all integrals are computed on the piecewise planar surface $\Gamma_{h,1}$ which is the zero level set of $\phi_{h,1}$ defined as in \eqref{eq:p1repoflsf}. See also \cite{OlRe18} where the mapping introduced in \cite{Le16} and its implementation is discussed.

\begin{figure}[!h]
	\centering
	\includegraphics{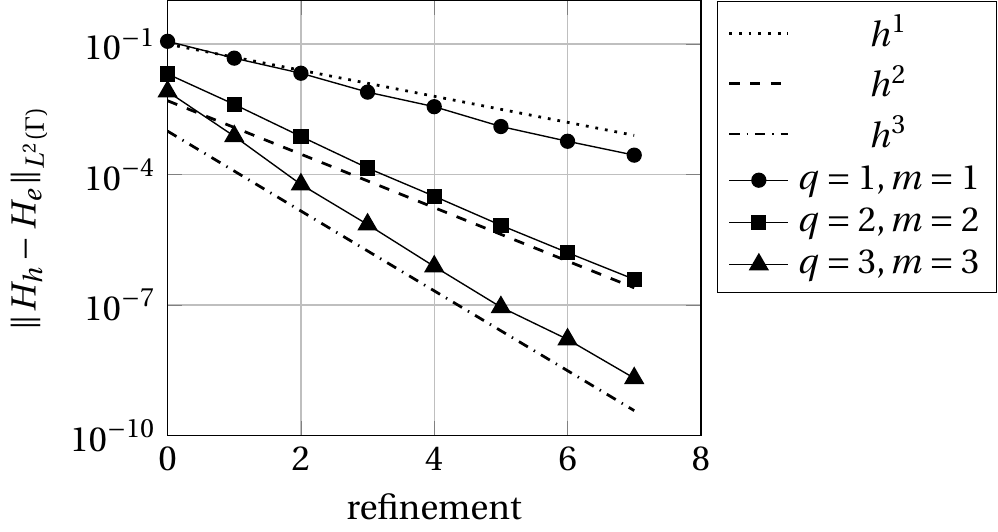}
	\caption{The convergence of the discrete mean curvature vector $\meancurv_h$ towards $H_e=\kappa \bfn$ in $L^2$ using the proposed stabilized finite element method. The level set function $\phi=\frac{x^2}{0.64}+y^2-0.25$. $m$ is the degree of the polynomials in the approximation space of $\meancurv_h$ and $q$ is the degree of the polynomials in the approximation space of $\phi_{h,q}$.\label{fig:curvatureError}}
\end{figure}

Given $\phi=\frac{x^2}{0.64}+y^2-0.25$ we compute the stabilized discrete mean curvature vector following the proposed strategy. The computational domain is taken to be $\Omega = [-1.2,  1.2]\times [-1.2, 1.2]$ and on this domain we generate a background mesh with an initial mesh size of $h = \frac{2.4}{10}$. At every refinement the mesh parameter $h$ is halved. We see in Fig. \ref{fig:curvatureError} that we obtain $q$th-order accurate approximations to $\kappa \bfn$ in $L^2$ with $m=q$.

In the numerical examples in the next section we use $q=2$ and $m=2$ and expect to have a second order accurate surface tension force. 

\section{Numerical examples}
We consider three numerical examples. In the examples in two space dimensions we use piecewise linear elements in time, i.e. $r=1$ while in three space dimensions we use piecewise constant functions in time, i.e. $r=0$. Note that using $r=0$ is equivalent to using the backward Euler method for the time discretization, see Remark 3.1. For the benchmark problems in two space dimensions we studied both the trapezoidal rule and the Simpson's quadrature rule in time. Both quadrature rules give similar results, therefore we only show the results using one of them. The results shown in the figures are with Simpson's quadrature rule. The stabilization parameters are $C_p=10^{-1}$ and $C_{\bfu}=10^{-2}$. For the surface tension force we use $q=2$, $m=2$, and $c_{F,j}=c_{\Gamma,j}=10^{-2}$ and expect to obtain a second order accurate approximation. 
We use $b=b^2$ in the discrete formulation \eqref{eq:weakformquad} since although the form $b^1$ in equation \eqref{eq:formb1} is mathematically equivalent to the form $b^2$ in equation \eqref{eq:formb2} we get lower spurious velocities with the form $b^2$, see Remark 5.1 in \cite{HaLaZa14}.

\subsection{Bubble in a pure straining flow}
We consider the example from \cite{example-Nasa}, section 5.2. We simulate the evolution of a bubble placed in a two dimensional slow viscous flow where inertial effects are negligible, $\bff=0$, but surface tension forces are important. Initially the bubble is a circle and \[
	\bfu(0,\bfx) = (Q x, -Qy),
\]
where $Q$ is the rate of shear. The bubble will evolve towards a stable steady state solution for $Q$ such that $0<Q<Q_c$, where $Q_c$ is the critical value of the rate of shear. The steady state solution in an infinite domain was first derived in \cite{strain-flow-exact-sol}.

For the computations we use a bounded domain $\Omega = [-L,L] \times [-L,L]$ with $L>0$ and we prescribe the following Dirichlet boundary conditions,
\[
	\bfu(t,\bfx) = (Q x, -Qy) \quad \forall \bfx \in \partial \Omega \ , \ \forall t > 0.
\]
The initial interface is a circle of radius $r_0 = 0.5$ centered in $(0,0)$. In order to measure the deformation of the bubble we use the deformation parameter 
\begin{equation}
D = \frac{R_{\max} - R_{\min}}{R_{\max} + R_{\min}},
\end{equation}
with $R_{\min}$ and $R_{\max}$ the minimum and the maximum distance of a point on the interface to the center of the bubble. The steady state solution for an unbounded domain can be explicitly computed through the following formula
\[
	Qa -2 I_0b = 0
\]
where $a$ and $b$ are real numbers describing the shape of the bubble at the steady state and $I_0$ is explicitly known and depends only on $a$ and $b$, see \cite{example-Nasa}.
\begin{figure}[!ht]
        \centering
         \subfloat
    	{
    	\scalebox{0.5}    	
	\centering	
	\includegraphics[scale=0.45]{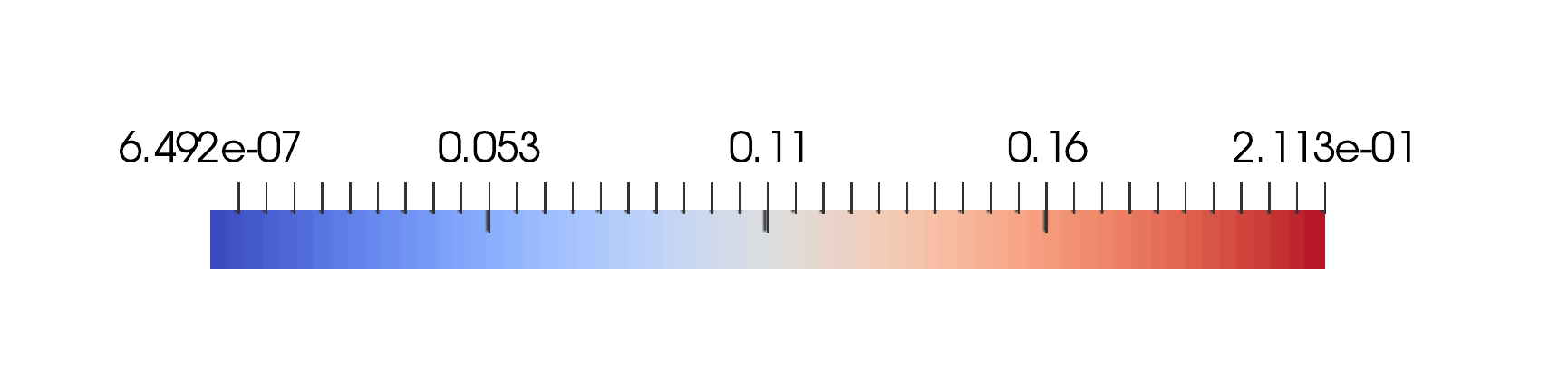}
        }
         \subfloat
    	{
    	\scalebox{0.5}    
	\centering	
        \includegraphics[scale=0.45]{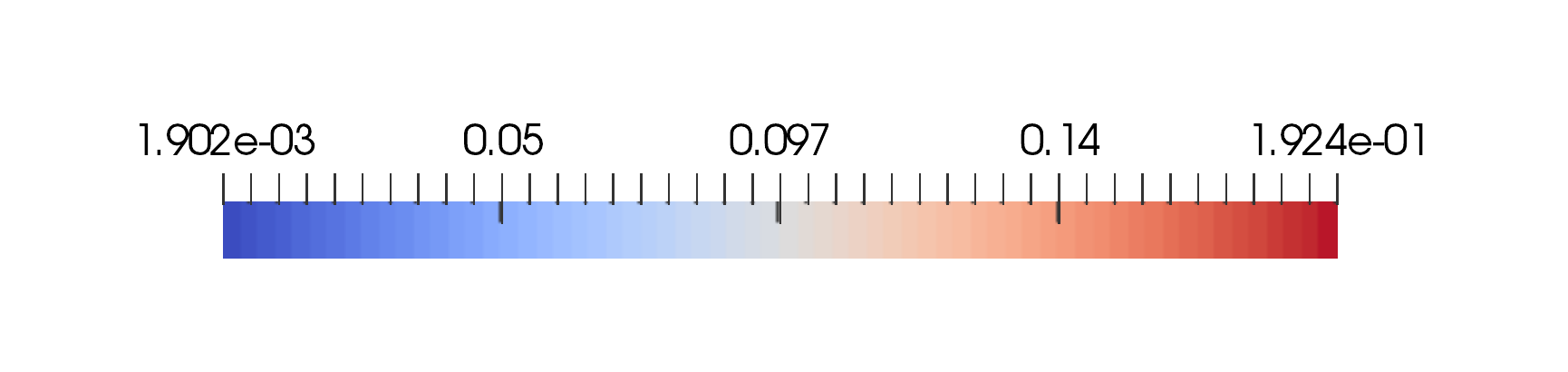}
        }\\[-0.8cm]
                \setcounter{subfigure}{0}

    \subfloat[$t=0$]
    	{
    	\scalebox{0.5}    	
	\centering	
	\includegraphics[scale=0.28]{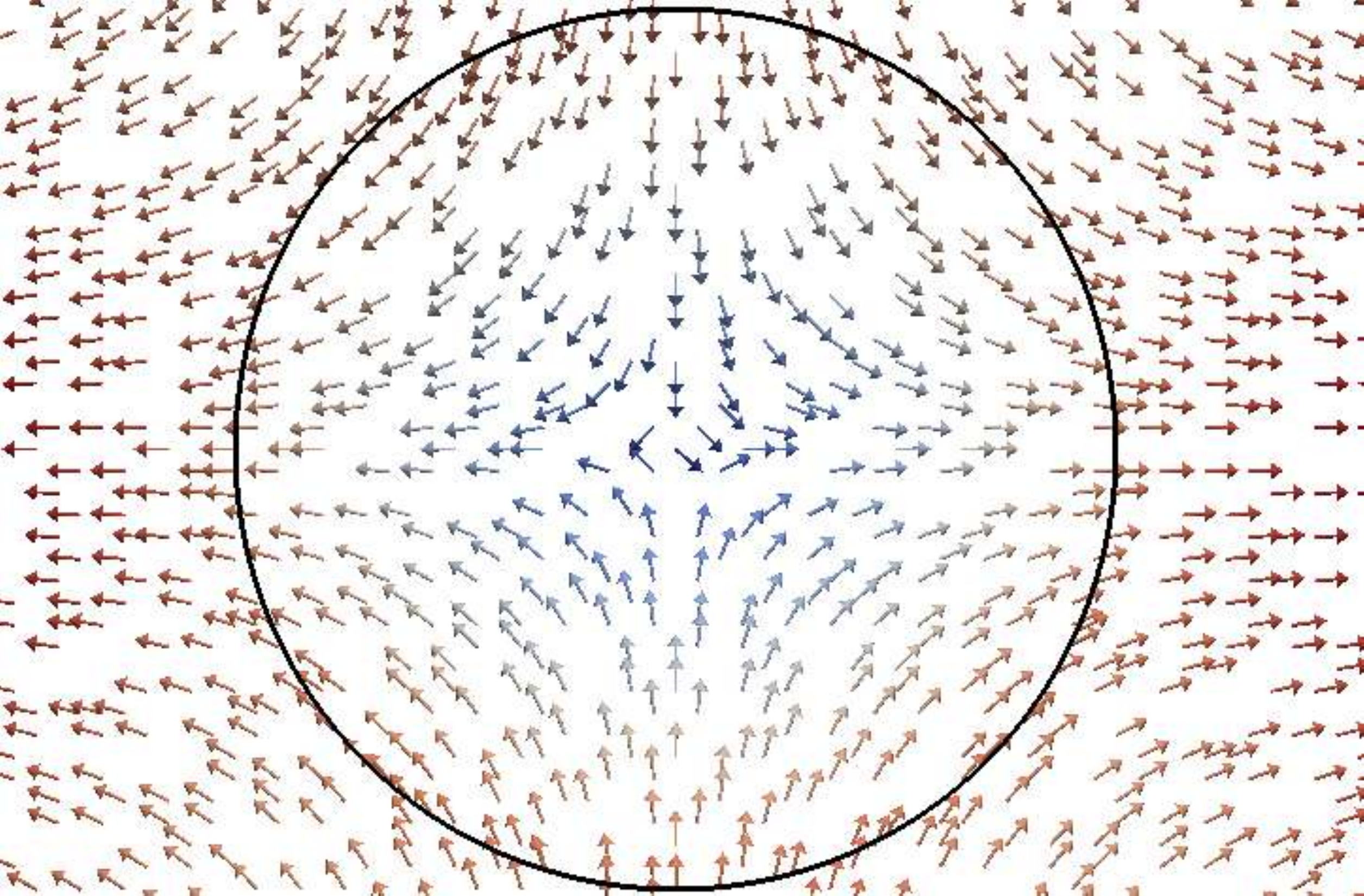}
        }
    \subfloat[$t=0.5$]
    	{
    	\scalebox{0.5}    	
	\centering	
        \includegraphics[scale=0.28]{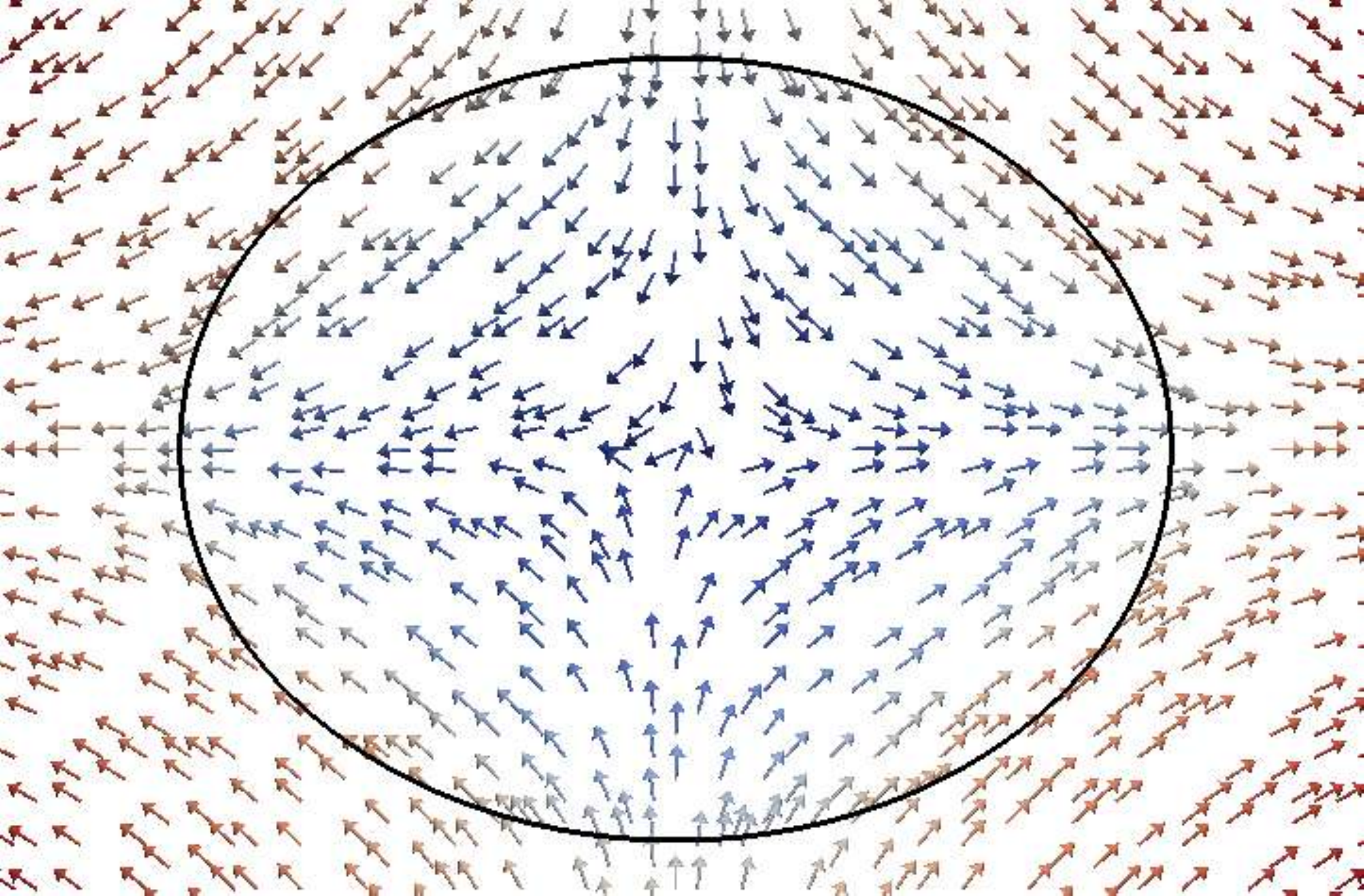}
        }
         \\[0.5cm]
         \subfloat
    	{
    	\scalebox{0.5}    	
	\centering	
	\includegraphics[scale=0.45]{colorBarDrop4HoriNoTitle}
        }
         \subfloat
    	{
    	\scalebox{0.5}    
	\centering	
        \includegraphics[scale=0.45]{colorBarDrop4HoriNoTitle}
        }\\[-0.8cm]
                        \setcounter{subfigure}{2}
    \subfloat[$t=1$]
    	{
        \scalebox{0.5}    	
	\centering	
	\includegraphics[scale=0.28]{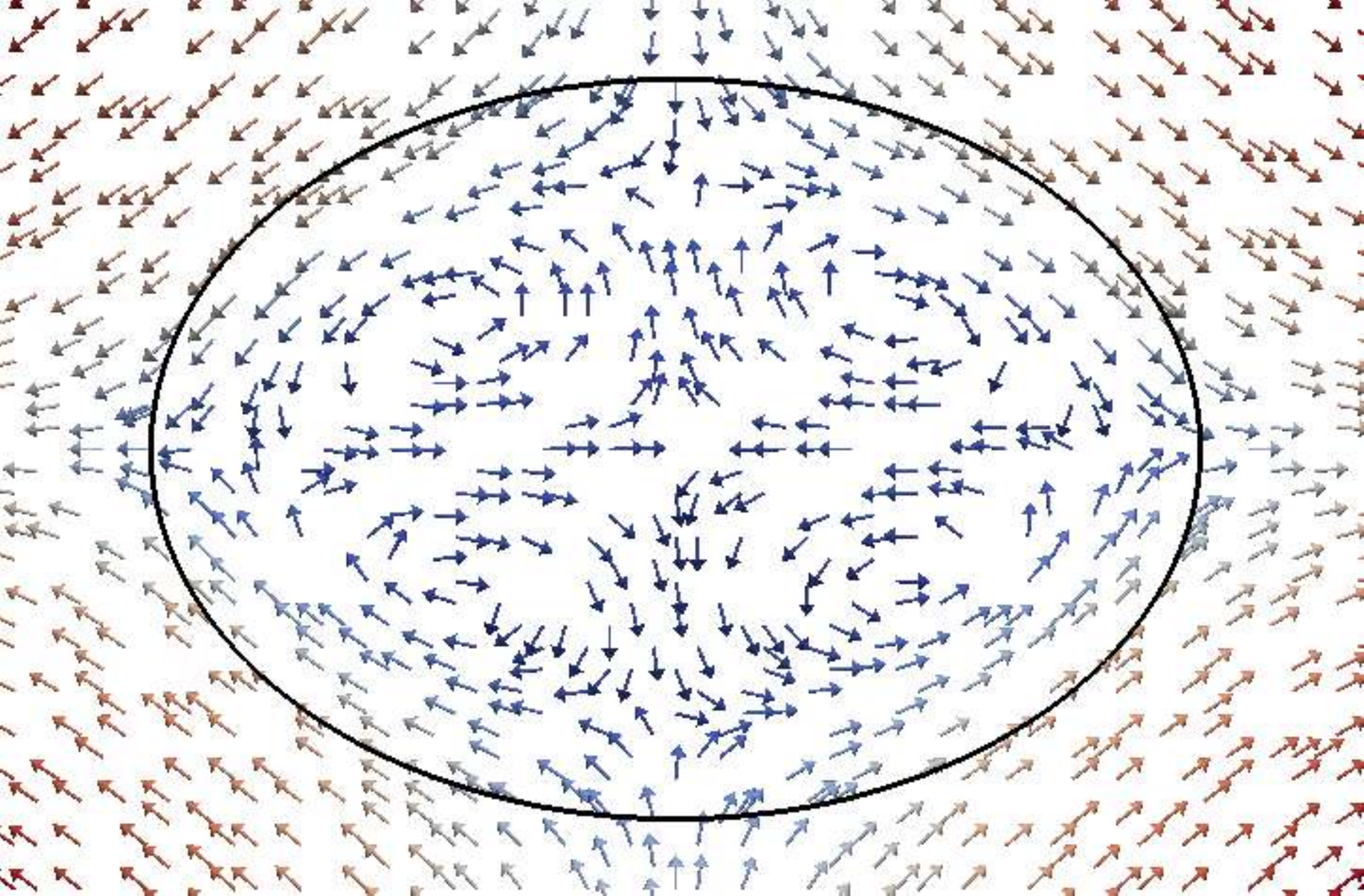}
        }   
         \subfloat[$t=1.6$]
    	{
        \scalebox{0.5}    	
	\centering	
	\includegraphics[scale=0.28]{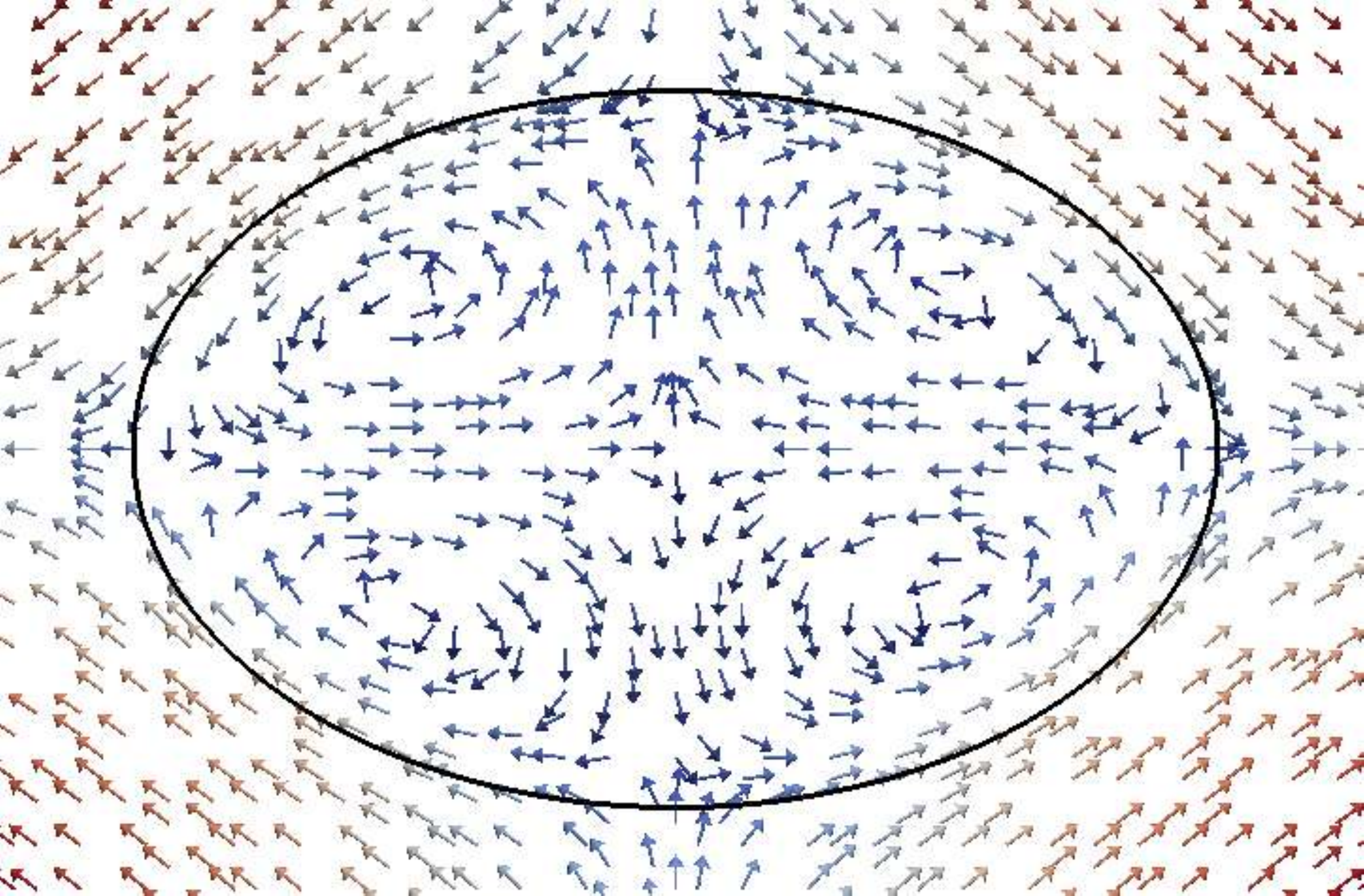}
        }   
\caption{Computed velocity fields at different times for a rate of shear $Q=0.2$ using a computational domain with $L=2$ and the proposed CutFEM. The mesh size is $h=1/40$ and the time step size is $\Delta t = h/4$.  The color bar shows the magnitude of the velocity field. \label{fig - strainFlow_drop1}}
\end{figure}

\begin{figure}[h]
	\centering  	
	\includegraphics{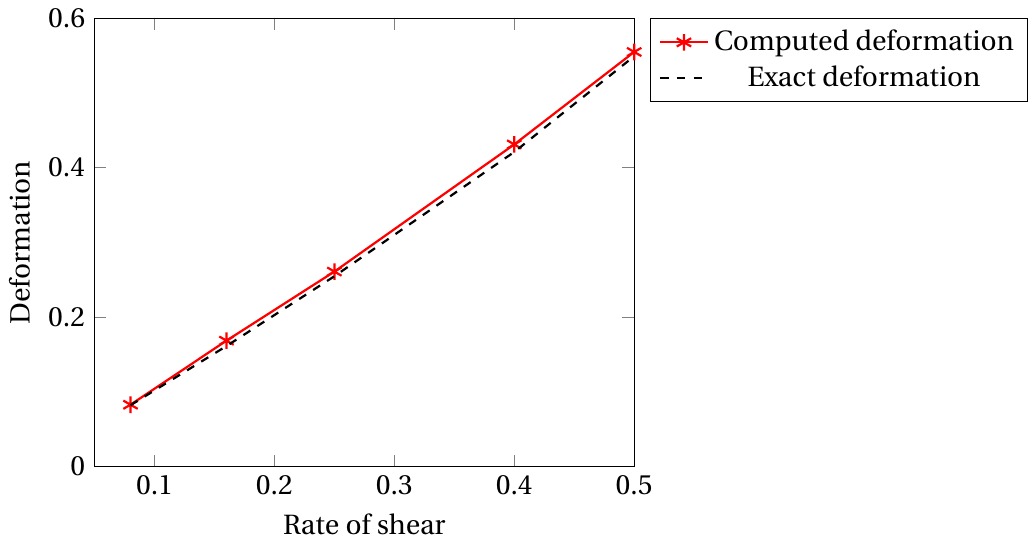}
	\caption{Exact (black dashed line) and the computed (red stars) deformation parameter $D$ as a function of the rate of shear $Q$. The computed deformation is shown for when the normal velocity is less than $1\times 10^{-5}$. The computation has been done using the proposed CutFEM in a computational domain with size $L=4$. The mesh size is $h=L/120$ and the time step size is $\Delta t = h / 4$. \label{fig:RateOfShearVsDeformation}}
\end{figure}

\begin{figure}[!h]
	\centering
	\includegraphics{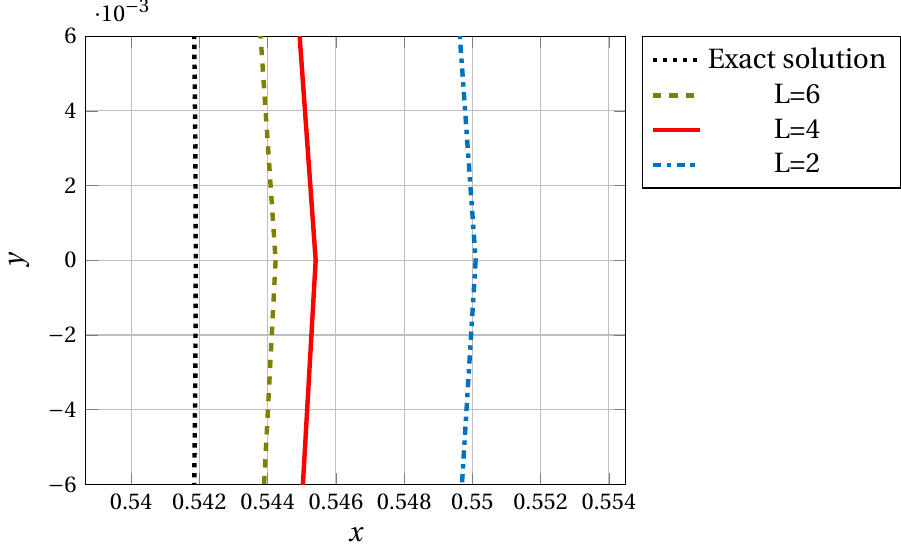}
    \caption{The exact solution in black and the computed solution using different sizes for the computational domain $\Omega = [-L,L] \times [-L,L]$. Blue: $L = 2$. Red: $L=4$. Green: $L=6$. The mesh size is $h=0.1$. 
\label{fig:stokes_interface_diffBox}}   	
\end{figure}

In Fig. \ref{fig - strainFlow_drop1} we show the computed bubble and the velocity field, using the proposed CutFEM, at different time instances and with $Q=0.2$. One can see that at the interface the velocity field becomes tangential to the interface and the shape of the bubble stops to change. For five different values of $Q$ we show the deformation of the bubble using the proposed CutFEM and the analytical solution in Fig. \ref{fig:RateOfShearVsDeformation}. The shown results are for when the normal velocity is less than $1\times 10^{-5}$. We see good agreement between our numerical results and the analytical solution. Since our computations are done on a bounded computational domain the obtained steady state solution is affected by the boundary. In Fig. \ref{fig:stokes_interface_diffBox}, we show the computed solution for different values of $L$ and we see that the solution comes closer to the analytical solution as the size of the computational domain increases.

\subsection{Rising bubble in 2D}\label{section - rising bubble 2D}
\begin{figure}[h]
	\centering
	\includegraphics{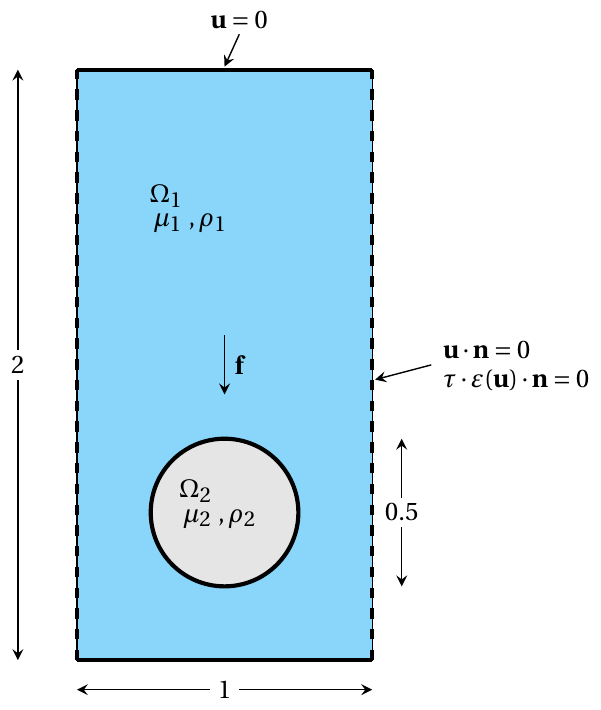}
	\caption{Initial configuration \label{fig:InitialConfiguration}}
\end{figure}
We now consider the benchmark test cases from \cite{Benchmark}.
A two-dimensional bubble rising in a liquid column due to gravity i.e., $\mathbf{f} = \rho (0,0.98)$. 
The computational domain in space is $[0,1] \times [0,2]$ and the bubble is initially a circle centered at $(0.5, 0.5)$ with radius $r_0 = 0.25$. The no-slip boundary condition, $\bfu = 0$, is imposed on the horizontal walls and the free slip condition, $\bfu \cdot \bfn = 0$ , $\tau \cdot 2\varepsilon(\bfu) \bfn = 0$, is imposed on the vertical walls. For an illustration of the initial configuration see Fig. \ref{fig:InitialConfiguration}. For the implementation of these boundary conditions in the proposed CutFEM see Remark \ref{rem:mixedBC} with $\mathbf{g}=0$.  
Two different test cases, with the physical parameters set as in Table \ref{table - set of parameters}, are considered.  
\begin{table}[h]
\centering
\begin{tabular}{ p{3cm} p{1cm} p{1cm} p{1cm} p{1cm} p{1cm}}
\hline 
Test case  & $\rho_1$ & $\rho_2$ & $\mu_1$ & $\mu_2$ & $\sigma$\\ 
\hline
1	& 1000	&	100	&	10	&	1	&	24.5 \\
2	& 1000	&	1	&	10	&	0.1	&	1.96 \\
\hline 
\end{tabular}
\captionof{table}{Parameters used in the two benchmark test cases. \label{table - set of parameters}}
\end{table}

To verify the accuracy of the proposed method we track the evolution of the bubble until time $T=3$, measure three benchmark quantities, defined below, and compare our results with results reported in \cite{Benchmark} by three other groups. See Table \ref{table - other methods} and  \cite{Benchmark} for more details about the methods used by the different groups. As in \cite{Benchmark} we use the following benchmark quantities:
\begin{itemize}
\item Center of mass 
\begin{align*}
\mathbf{X}_c = (x^1_c, x^2_c) = \frac{\int_{\Omega_2} \bfx \diff x}{\int_{\Omega_2} 1 \diff x}.
\end{align*}
Here, the second component $x^2_c$ is of interest.
\item Circularity 
\begin{align*}
c = \frac{P_a}{P_b},
\end{align*}
where $P_a$ is the perimeter of the circle which has an area equal to that of the bubble with perimeter $P_b$.
\item Rise velocity 
\begin{align*}
\mathbf{U}_c = (u^1_c, u^2_c) = \frac{\int_{\Omega_2} \bfu \diff x}{\int_{\Omega_2} 1 \diff x}.
\end{align*}
Here, we use the velocity component $u^2_c$ which is in the direction opposite to the gravitational vector $\bff$.
\end{itemize}
In order to measure the error in those quantities, the following relative error norms are used
\begin{align}
l_1 \text{ norm - } \norm{e}_1 & = 
\frac{\sum_{t=1}^{NTS} | \omega_{t,\text{ref}} - \omega_t| }
     {\sum_{t=1}^{NTS} | \omega_{t,\text{ref}}| },  \label{eq:l1norm} \\
l_2 \text{ norm - } \norm{e}_2 & = \left(
\frac{\sum_{t=1}^{NTS} | \omega_{t,\text{ref}} - \omega_t|^2 }
     {\sum_{t=1}^{NTS} | \omega_{t,\text{ref}}|^2 } \right)^{\frac{1}{2}},  \\ 
l_{\infty} \text{ norm - } \norm{e}_\infty & = 
\frac{\max_t | \omega_{t,\text{ref}} - \omega_t| }
     {\max_t | \omega_{t,\text{ref}}| },   \label{eq:linfnorm} 
\end{align}
where $\omega_t$ denote the benchmark quantity at time instance $t$. The reference solution  $\omega_{t,\text{ref}}$ is the solution at time instance $t$ computed on the finest grid and NTS is the number of time steps.
The time step size was chosen equal to $\Delta t = h/4$.

\begin{table}[h]
\centering
\begin{tabular}{ p{2.5cm} p{5cm} p{3.2cm} p{4cm}}
\hline 
Group  & Method  & Interface & Time discretization \\ 
\hline
1 : TP2D		& Unfitted FEM, $\tilde{\mathbb{Q}}_1\mathbb{Q}_0$	&	Level-set $\mathbb{Q}_1$	& Fractional step $\Theta$ scheme	\\
2 : FreeLIFE	& Unfitted FEM, $\Po_1-\text{iso}-\Po_2 / \Po_1$ & Level-set $\Po_ 1$ & BDF2\\
3 : MooNMD 		& Fitted FEM, $\Po_2$ enriched with cubic polynomials- discontinuous $\Po_1$ & Lagrangian markers 
& Second order fractional step $\Theta$ scheme\cite{Ran04}\\		
\hline 
\end{tabular}
\captionof{table}{The different groups from \cite{Benchmark} and the different computational techniques that are used. \label{table - other methods}}
\end{table}

\subsubsection{Benchmark test case 1}
\begin{figure}[h]
\centering
    \subfloat
    	{
	\includegraphics[scale=0.9]{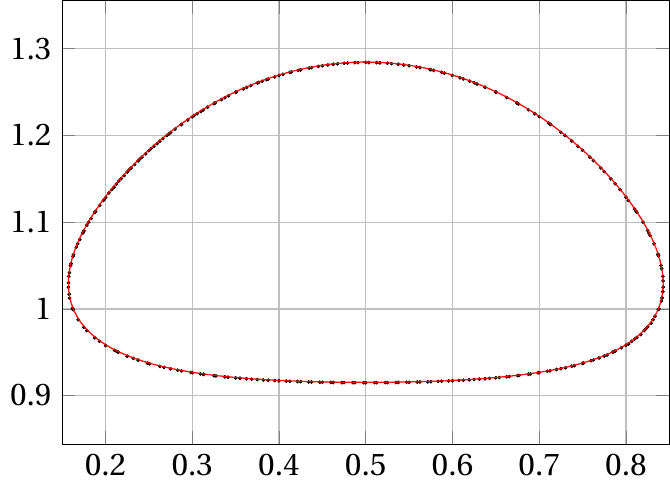}
        }  
      \hspace{1.5cm}
    \subfloat
    	{
	\includegraphics[scale=0.85]{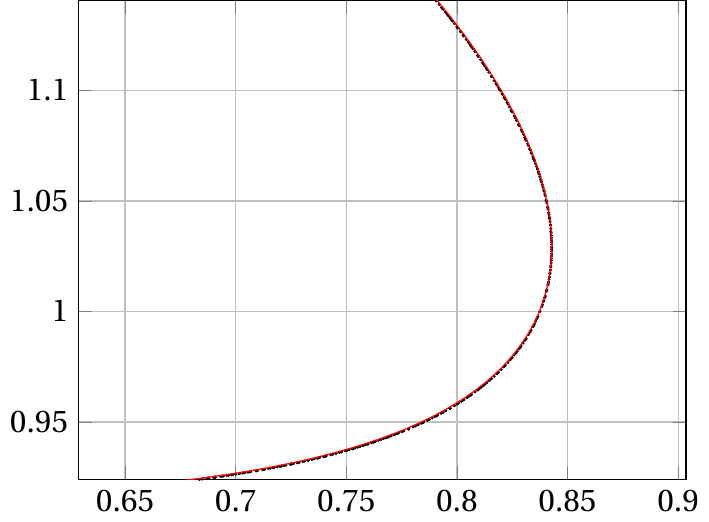}
        }       
        \caption{Left panel:  The shape of the bubble at time $t=3$ obtained using the proposed CutFEM. The shape of the bubble on the course mesh, $h=\frac{1}{40}$, (black dots) is compared with the shape obtained on the finest mesh, $h=\frac{1}{160}$, (red line). Right panel:  comparison of the shape of the bubble computed on our finest mesh (red line) with the shape obtained by group 1 in \cite{Benchmark} for the same mesh size (black dots).
\label{interface1_plot}
        }   	
\end{figure}

We first look at the bubble shape obtained at the final time $t=3$. In Fig.~\ref{interface1_plot} (left panel) we see that the solution on the coarse mesh, $h=1/40$, and the solution on the finest mesh, $h=1/160$, are not distinguishable. Moreover, we show a close up of the bubble shape on our finest grid and the solution of group 1 on the same grid. Again no difference between the two shapes are visible.

\begin{figure}[h]
\centering
    \subfloat
    	{
	\includegraphics[scale=0.8]{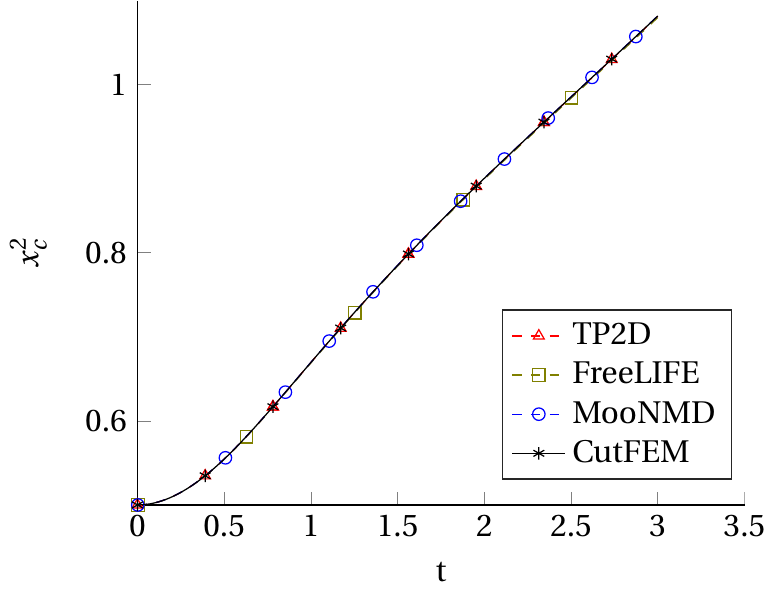}
        }  
    \subfloat
    	{
	\includegraphics[scale=0.8]{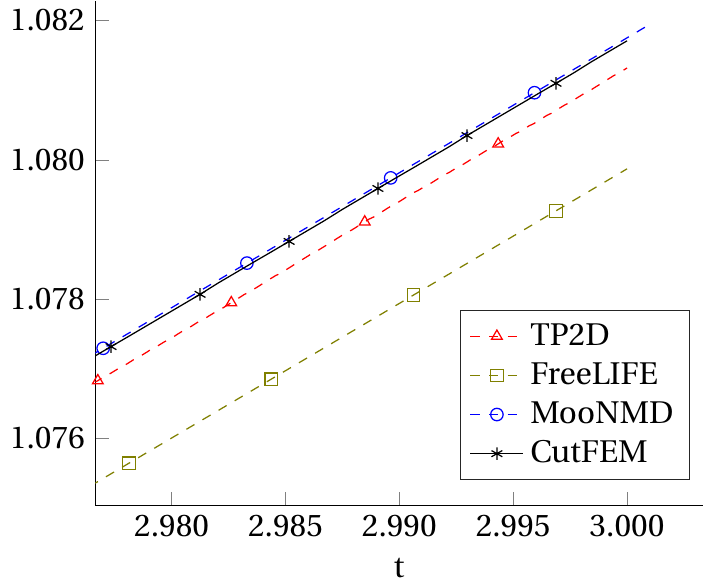}
        }       
        \caption{Center of mass as a function of time. Left panel: the center of mass obtained by the proposed CutFEM compared with the results of the three groups in \cite{Benchmark}. Right panel: close-up of the center of mass around the final time.\label{fig:centerOfMass1}
        }   	
\end{figure}

\begin{figure}[!h]
\centering
    \subfloat
    	{
	\includegraphics[scale=0.8]{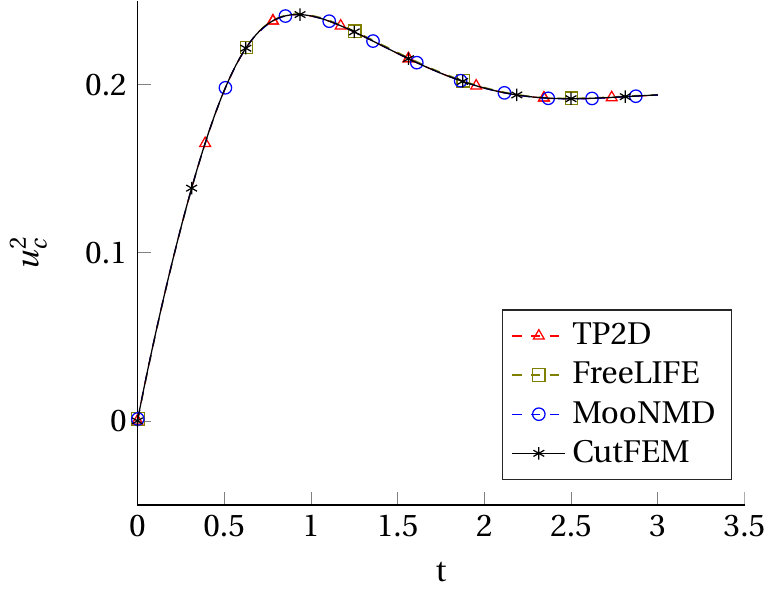}
         }
    \subfloat
    	{
	\includegraphics[scale=0.8]{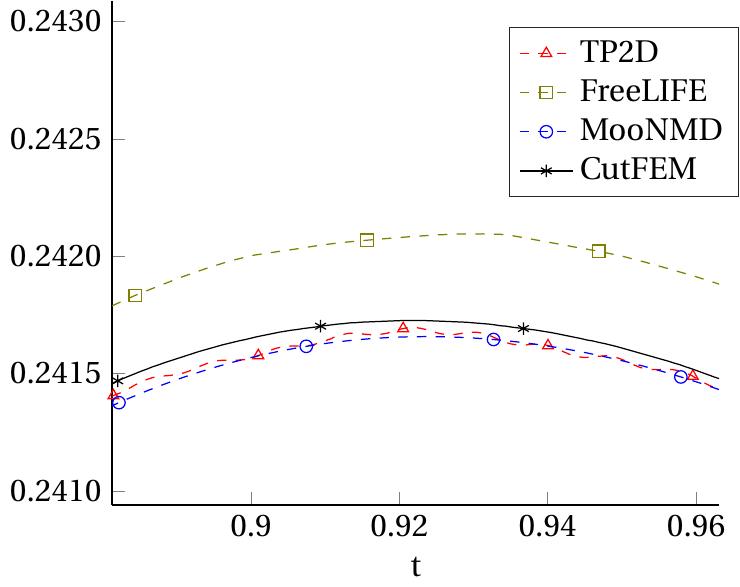}
        }
        \caption{Rise velocity as a function of time. Left panel: rise velocity obtained by the proposed CutFEM compared with the results of the three groups in \cite{Benchmark}. Right panel: close-up of the rise velocity where the rise velocity is maximal. \label{fig:riseVelocity1}
        }   	
\end{figure}

\begin{figure}[!h]
\centering
    \subfloat
    	{
        \scalebox{0.5}    	
	\centering	
	\includegraphics[scale=0.8]{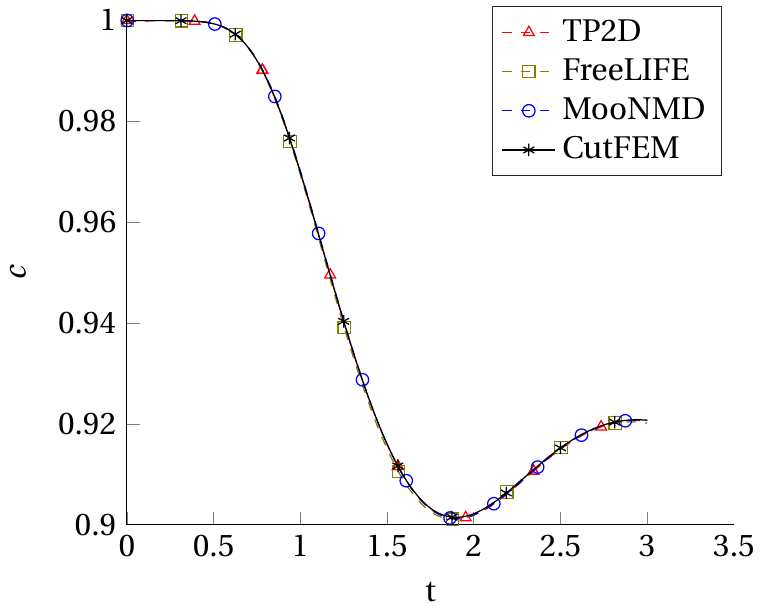}
        }  
        \addtocounter{subfigure}{-1}
        \\
    \subfloat[]
    	{
    	\scalebox{0.3}    	
	\centering	
	\includegraphics{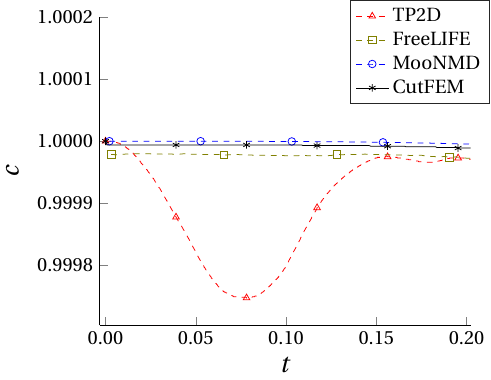}
        }
    \subfloat[]
    	{
    	\scalebox{0.3}    	
	\centering	
	\includegraphics{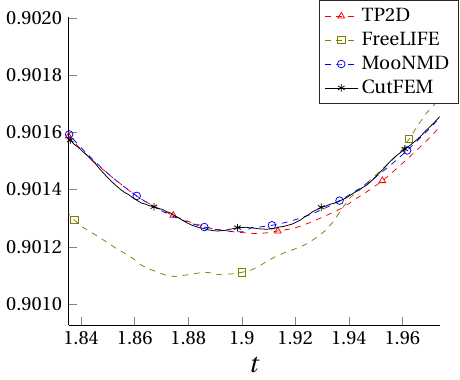}
        }
    \subfloat[]
    	{
        \scalebox{0.3}    	
	\centering	
	\includegraphics{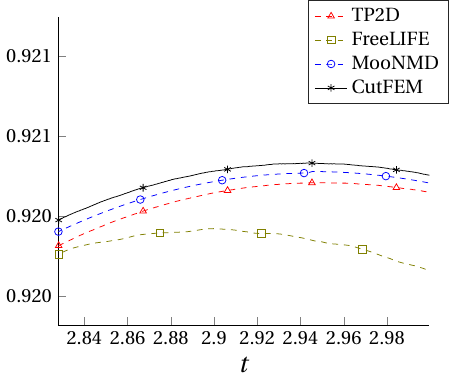}
        }                   	
        \caption{Circularity as a function of time. Bellow: Close up of the circularity, (a) at the initial time, (b) where the deformation of the bubble is maximal, and (c) around the final time.   \label{fig:circularity1_zoom1-2}}   	

\end{figure}

\begin{figure}[!h]
\centering
	\includegraphics[scale=1]{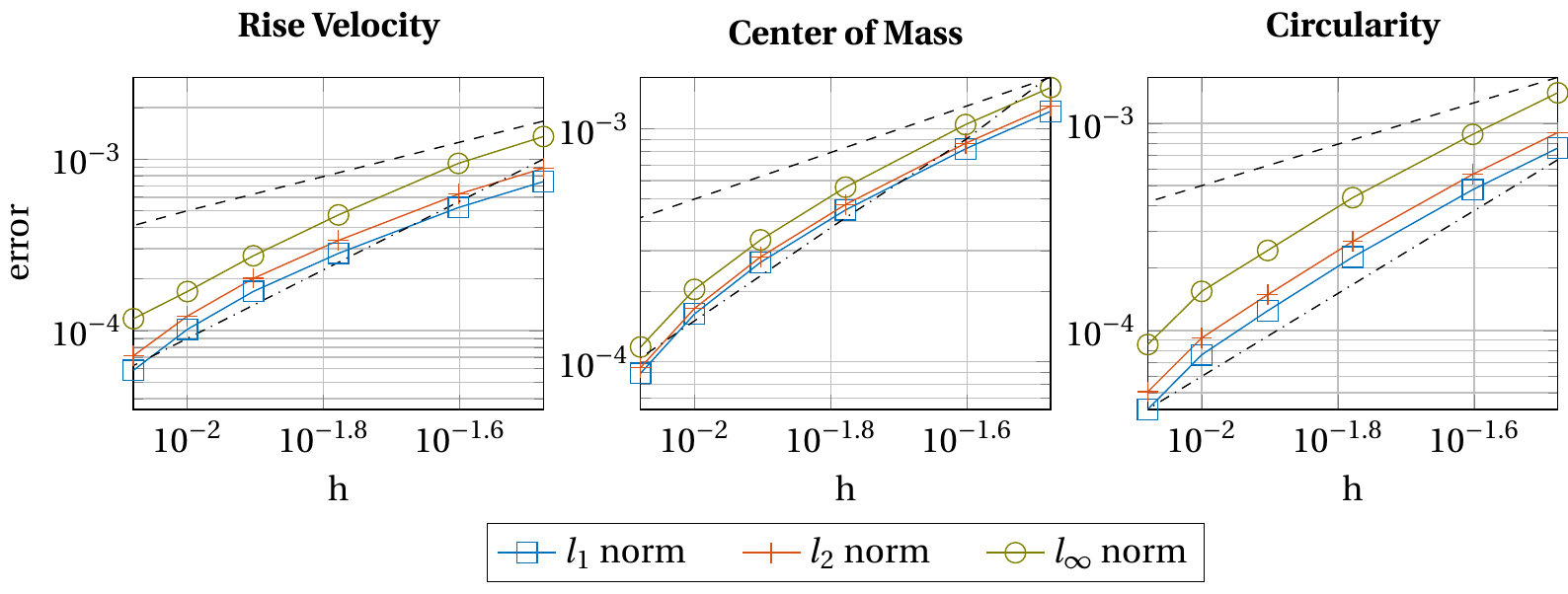}            
        \caption{Convergence of the benchmark quantities, the rise velocity, center of mass, and circularity, in the norms, $l_\infty$ (o), $l_2$ (+), and $l_1$ ($\square$). The dashed line is proportional to $h$ and the dash-dotted line is proportional to $h^2$. \label{convergence_Benchmark1}}  
\end{figure}

\noindent In Fig.~\ref{fig:centerOfMass1},~\ref{fig:riseVelocity1}, and~\ref{fig:circularity1_zoom1-2} we show the evolution of the three benchmark quantities. We observe good agreement with the groups from~\cite{Benchmark}. When zooming in on the different parts where we see some differences, we see that our results is closest to MooNMD (group 3). Note that in MooNMD the mesh is conformed to the evolving interface and thus a re-meshing process is used. 
\noindent
In Fig. \ref{convergence_Benchmark1} we show the convergence of the benchmark quantities in the different norms defined in \eqref{eq:l1norm}-\eqref{eq:linfnorm}. We see that for the coarse mesh the convergence order is slower than 2 but the convergence order increases when the mesh is refined and one can see that it is around 2 in all three norms. 

\begin{figure}[!h]
\centering
\includegraphics[scale=0.4]{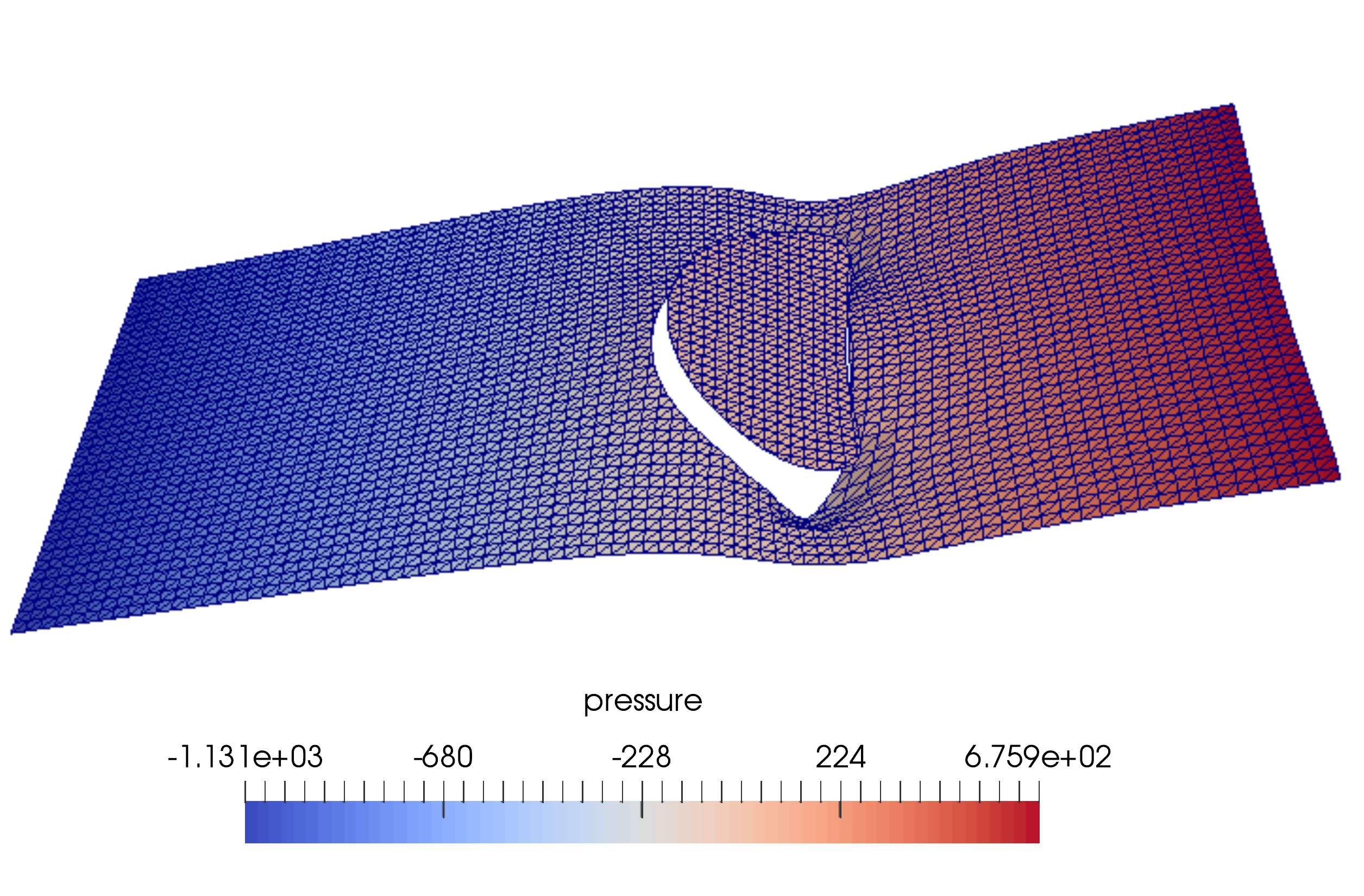}
        \caption{The discontinuous pressure at $t=1.5$ obtained using the proposed CutFEM. The mesh size is $h=1/40$ and the time step size is $\Delta t = h / 4$.\label{jumpKink_Benchmark1}}   	
\end{figure}

Finally, we show the discontinuous pressure at time $t=1.5$ on the course mesh, i.e. $h=1/40$, in Fig. \ref{jumpKink_Benchmark1}. We see that the proposed CutFEM can capture discontinuities without aligning the mesh to the interface.   
 
\subsubsection{Benchmark test case 2}
\begin{figure}[!h]
\centering
    \subfloat[$t=0$]
    	{
    	\scalebox{0.25}    	
		\centering	
		\includegraphics[scale=0.13]{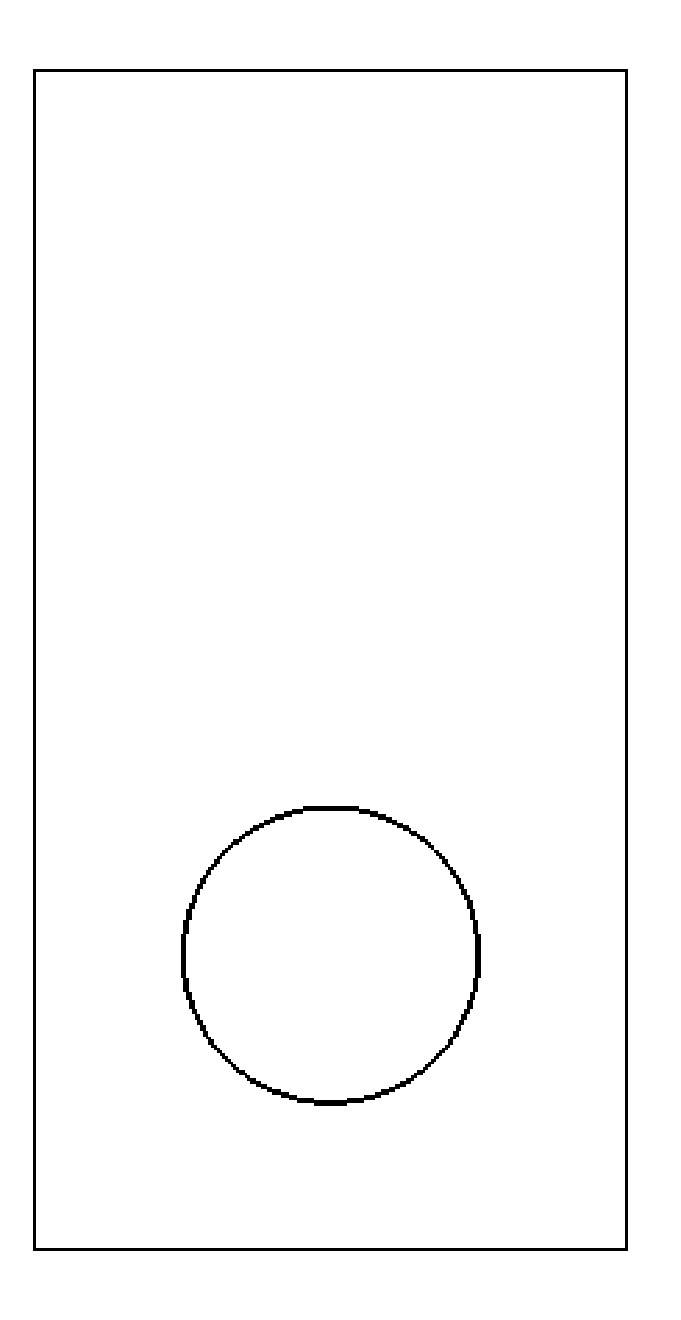}
        }
    \subfloat[$t=1$]
    	{
    	\scalebox{0.25}    	
		\centering	
		\includegraphics[scale=0.13]{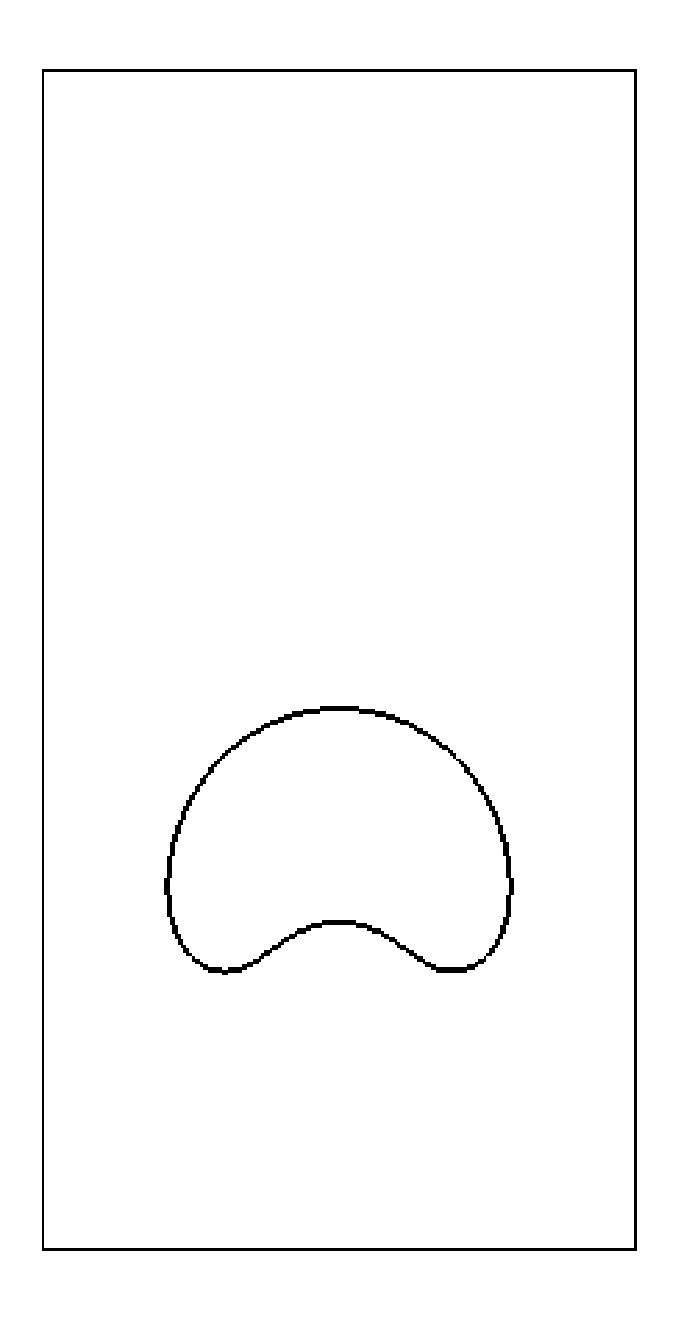}
        }
    \subfloat[$t=2$]
    	{
        \scalebox{0.25}    	
		\centering	
		\includegraphics[scale=0.13]{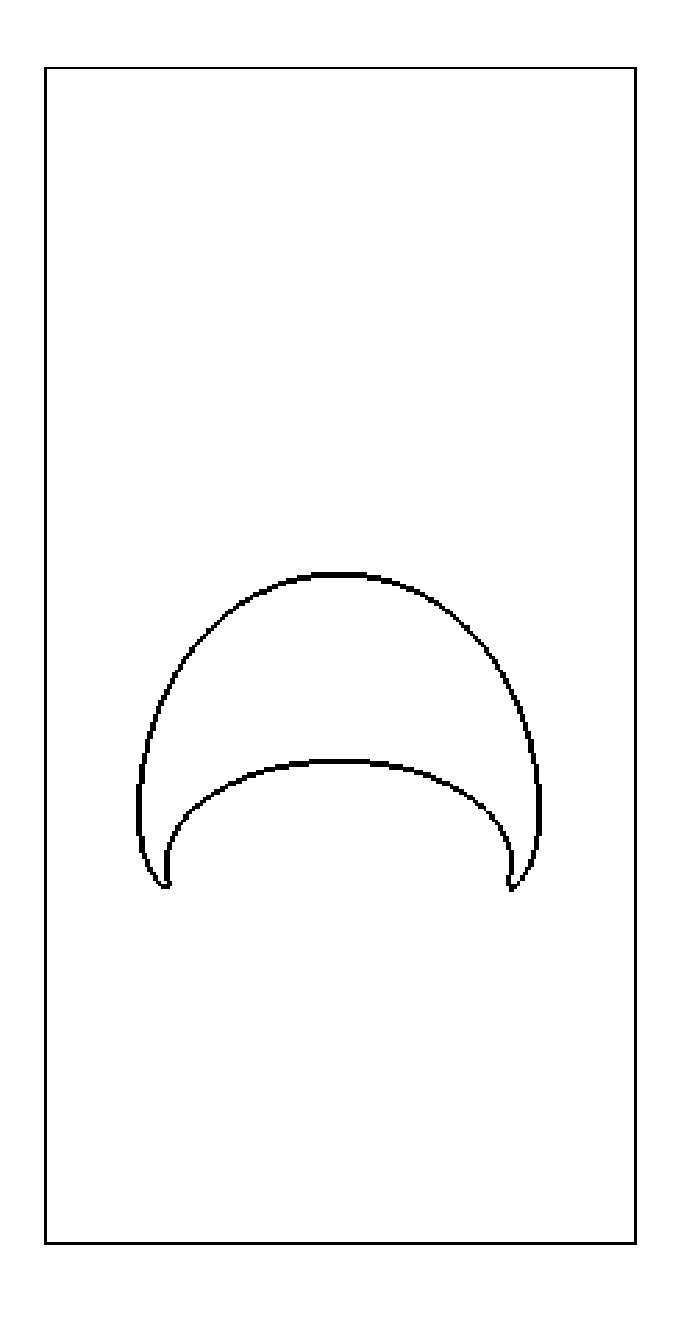}
        }    
    \subfloat[$t=3$]
    	{
        \scalebox{0.25}    	
		\centering	
		\includegraphics[scale=0.13]{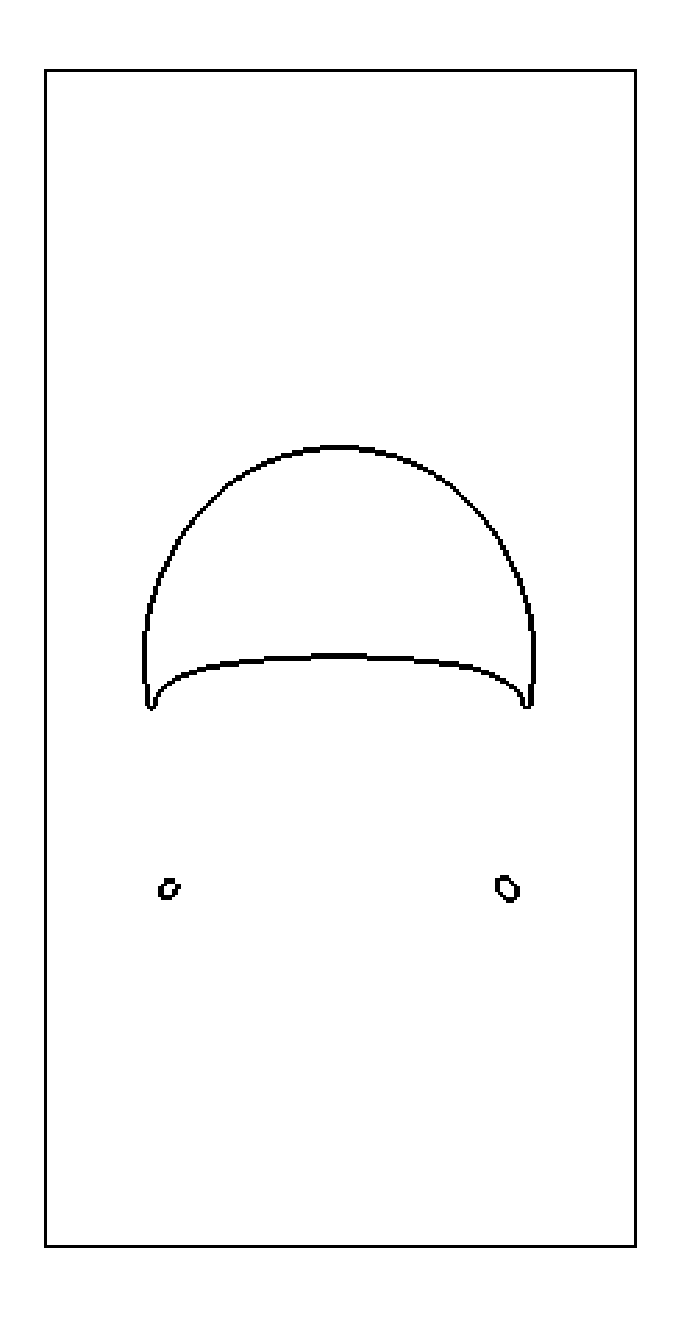}
        }                 	
        \caption{Bubble at different time instances. The mesh size is $h = 1/80$ and the time step size is $\Delta t = h / 4$. \label{fig:multiple shape}}   	
\end{figure}

\begin{figure}[!h]
	\centering
	\includegraphics{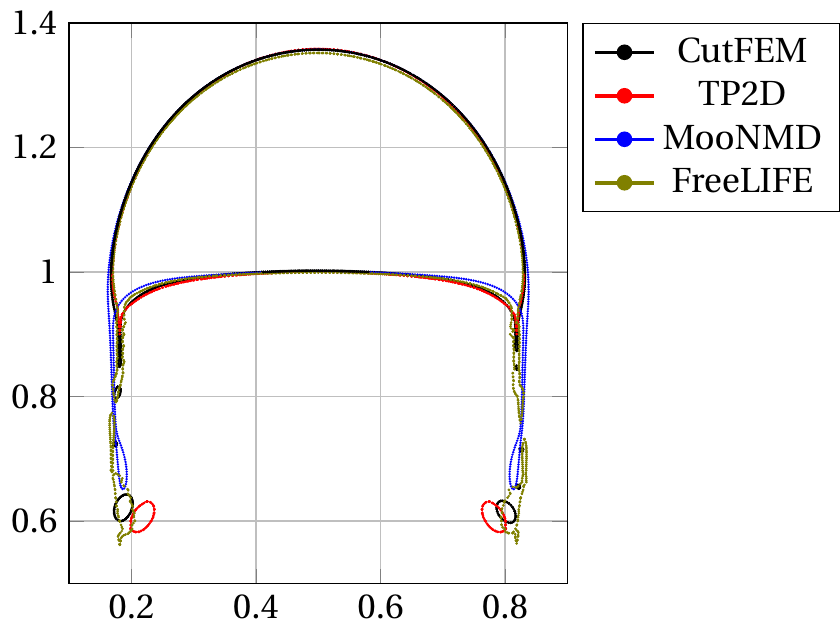}
	\caption{ Comparison of the shape of the bubble obtained, on the finest meshes, by CutFEM (black), TP2D (red), MooNMD (blue), and FreeLIFE (green).  \label{fig:interface2_plot}}
\end{figure}

In this test case, the low surface tension causes the bubble to deform, and the development leads to filaments and/or breaks up. 
One can see this behavior in Fig. \ref{fig:multiple shape} and \ref{fig:interface2_plot}. 
\begin{figure}[!h]
	\centering
	\includegraphics[scale = 1]{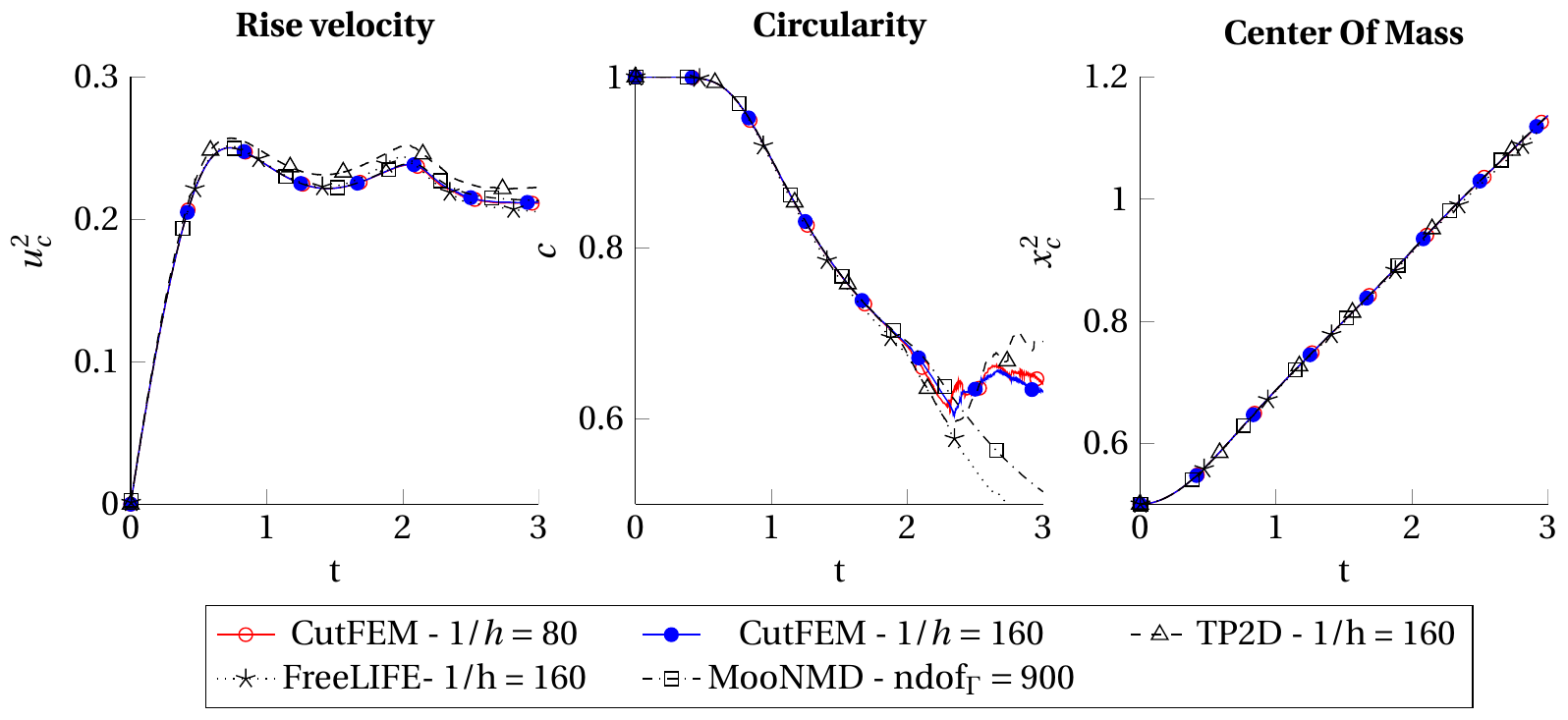}
	\caption{Center of mass, rise velocity, and circularity as a function of time. Results obtained with the proposed CutFEM are compared with the results of the three groups in \cite{Benchmark}. \label{fig:benchmark2}}
\end{figure}
The evolution of the different benchmark quantities are shown in Fig. \ref{fig:benchmark2}. We observe that the rise velocity obtained by TP2D is higher than what we and the other groups obtain. However, this difference vanishes when the mesh is refined and for $1/h=640$ group 1 also obtains similar rise velocity curves. Before the point of break up the computed circularity obtained by the different methods are similar. The proposed method as well as the method by group 1 end up with small satellites and break up after a while. 
We believe that the break up is caused by the numerical method we use for representing and evolving the interface and is due to the difficulty to correctly resolve the thin long filaments.

\subsection{Rising bubble in 3D}
In this three-dimensional example we simulate a gas bubble rising in a liquid. We consider different cases where we vary the fluid viscosities and the surface tension coefficient, see Table \ref{table - set of parameter 3D}.
\begin{table}[h]
\centering
\begin{tabular}{ p{3cm} p{1cm} p{1cm} p{1cm} p{1cm} p{1cm}}
\hline 
Shape  & $\rho_l$ & $\rho_g$ & $\mu_l$ & $\mu_g$ & $\sigma$\\ 
\hline
Spherical	& 1000	&	10	&	62  &	0.62	&	245  \\
Ellipsoidal	& 1000	&	10	&	35	&	0.35	&	24.5 \\
Skirted		& 1000	&	10	&	11	&	0.11	&	2.45 \\
\hline 
\end{tabular}
\captionof{table}{Parameters of the different cases. \label{table - set of parameter 3D}}
\end{table}

\begin{figure}[h!]
\centering
    \subfloat[Spherical shape]
    	{
    	\scalebox{0.3}    	
		\centering	
		\includegraphics[scale=0.22]{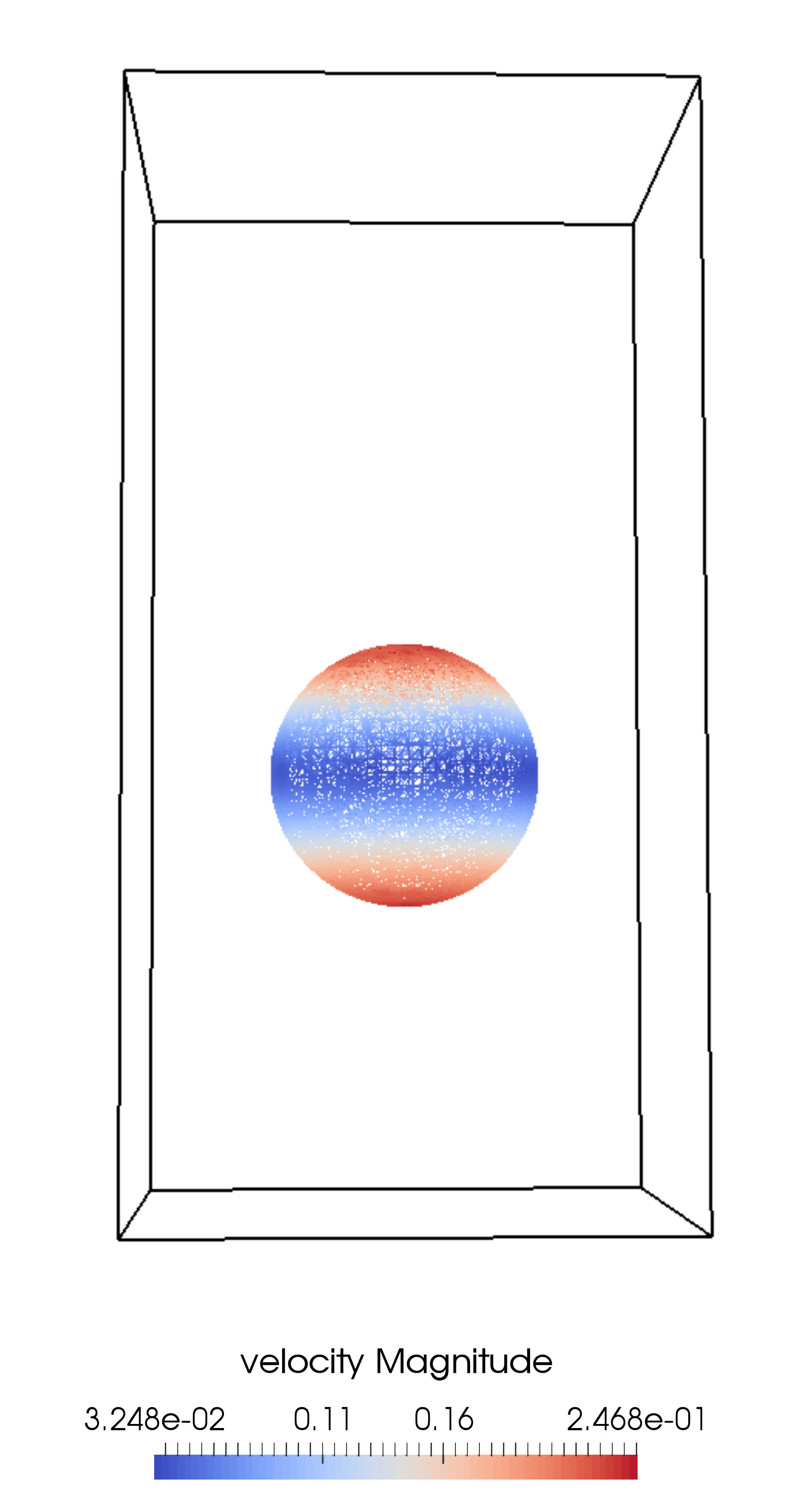}
        }
    \subfloat[Ellipsoidal shape]
    	{
    	\scalebox{0.3}    	
		\centering	
		\includegraphics[scale=0.22]{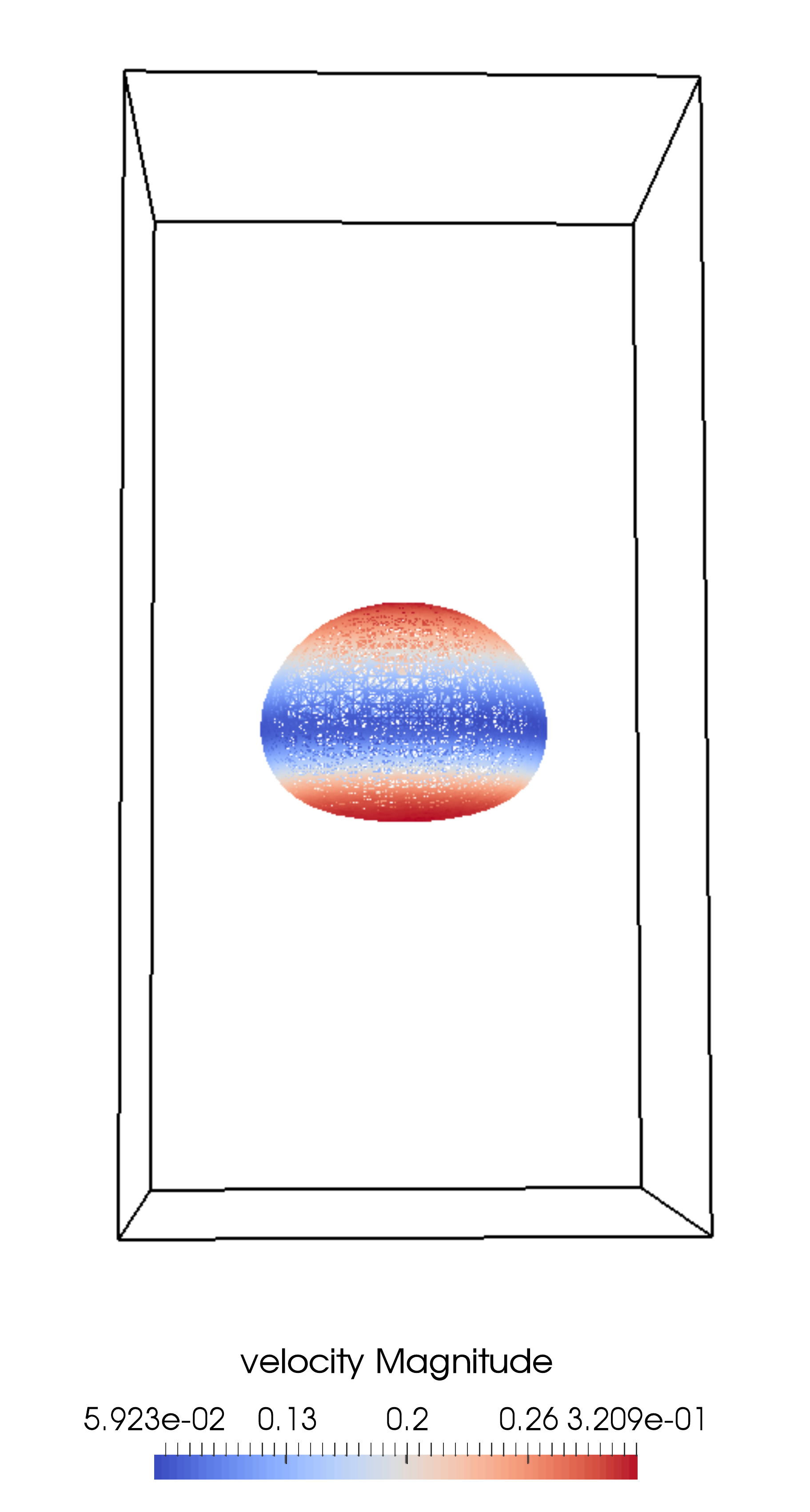}
        }
    \subfloat[Skirted shape]
    	{
        \scalebox{0.3}    	
		\centering	
		\includegraphics[scale=0.22]{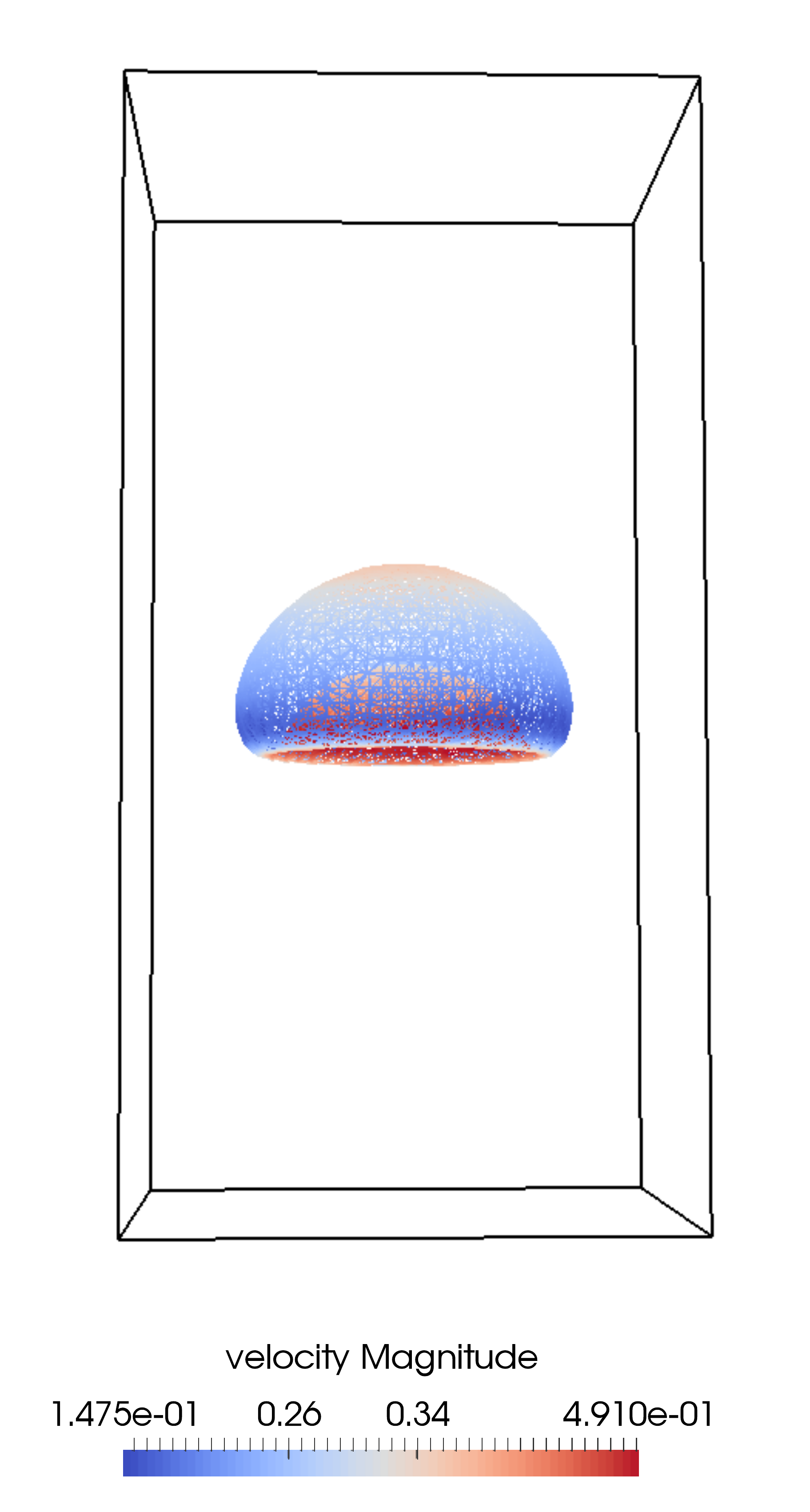}
        }                   	
        \caption{ Different shapes obtained at time $t=1.5$ using the different parameters in Table \ref{table - set of parameter 3D}  \label{rising_drop_3D}}   
\end{figure}

The computational domain in space is $[-0.5, \ 1.5] \times [-0.5, \ 1.5] \times [0, \ 2]$, the interface is initially a sphere of radius $r_0 =0.25$ centered in $(0.5, 0.5, 0.5)$. We use no-slip boundary condition on the horizontal walls and the free slip boundary condition is applied on vertical walls.\\
All the computations have been done on a grid of $40 \times 40 \times 80$ points and using a time step size $\Delta t = 2.5 \times 10^{-2}$. In Fig.~\ref{rising_drop_3D} we show the shape of the gas bubble at time $t=1.5$ for the different parameters in Table~\ref{table - set of parameter 3D}. We track the evolution of the bubble until time $T=3$ and evaluate the final shape of the bubble obtained by the proposed CutFEM with the diagram of Grace \cite{Grace}, which classifies the shapes of bubbles in quiescent viscous liquids. The governing dimensionless numbers are the Morton number (M), the Eotvos number (Eo) and the Reynolds number (Re) given by
\begin{align*}
M  &= \frac{g \mu_l^4 }{\rho_l \sigma^3}, &
Eo &= \frac{g  \Delta \rho d_e^2}{\sigma}, &
Re &= \frac{\rho_l v_\infty d_b}{\mu_l}, 
\end{align*}
where the subscript $l$ is used for the liquid parameters, $d_e$ is the diameter of a sphere with the same volume as the bubble, $d_b$ is the diameter of the bubble, and $v_\infty$ is the terminal rise velocity defined as in Section \ref{section - rising bubble 2D}. 
We compute the terminal rise velocity and thus the Reynolds number and present our results in Table \ref{table - result bubble 3D}.  We compare our results with the Reynolds number measured in \cite{Grace} but also with two other simulations, the level set/ghost fluid method in \cite{Benchmark-3D} and the volume of fluid method in \cite{Benchmark-3D_2}. We observe that even on a coarse mesh and with a large time step size, we obtain small differences between our computed Reynolds numbers and the Reynolds numbers measured by Grace. We use for example $\Delta t=\frac{\Delta x}{2}$ while in \cite{Benchmark-3D}  $\Delta t/\Delta x$ is $1.6 \cdot 10^{-4}$ in the spherical and ellipsoidal case and $1.6 \cdot 10^{-3}$ in the skirted case. 

\begin{table}[h!]
\centering
\begin{tabular}{ p{2cm} p{1.2cm} p{1cm} p{1.5cm} p{1.5cm} p{1.5cm} p{1cm}}
\hline 
Shape  & $M$ & $Eo$ &$Re_\text{Grace}$ & $Re_\text{CutFEM}$   & $Re_1 $ & $Re_2 $ \\ 
\hline
Spherical	& $1\times10^{-3}$	&	$1$		&	$1.7$	&	$1.77$ 	& 	$1.73$  &  $1.6$	  	\\
Ellipsoidal	& $0.1$				&	$10$	&	$4.6$	&	$4.63$ 	& 	$4.57$  &  $4.3$ 	\\
Skirted		& $1$				&	$100$	&	$20.0$	&	$19.8$  &	$19.21$ &  $18$ 	\\
\hline 
\end{tabular}
\captionof{table}{The dimensionless numbers for the different cases and comparison of the terminal computed Reynolds number from CutFEM ($Re_\text{CutFEM}$), the terminal Reynolds number computed in \cite{Benchmark-3D} ($Re_1$),  the terminal Reynolds number computed in \cite{Benchmark-3D_2} ($Re_2$), and the Reynolds number measured by Grace \cite{Grace} ($Re_\text{Grace}$). \label{table - result bubble 3D}}
\end{table}

\section{Conclusion}
We have presented a space-time cut finite element method which is able to accurately capture both the strong discontinuity in the pressure and the weak discontinuity in the velocity field across moving interfaces separating immiscible fluids, without conforming the mesh to these interfaces. We have also proposed a new  high order accurate finite element method for computing the mean curvature vector and consequently the surface tension force. Numerical experiments show that by stabilizing the $L^2$ projection we can compute a stabilized mean curvature vector based on the Laplace-Beltrami operator with improved accuracy compared to not stabilizing. The presented space-time method has also a convenient implementation as it does not reconstruct the space-time domain but rather directly uses quadrature rules to approximate the space-time integrals in the variational formulation.  The time discretization in the proposed method is closely related to implicit finite difference methods and we showed that the backward Euler method can be obtained by using piecewise constant functions in time.  In this paper we have presented a method which yields a second order accurate velocity approximation. However, are aim has also been to present a strategy which can be extended to yield higher order approximations if the regularity of the problem allows. 

In this work we combined the presented CutFEM with a level set method for the representation and evolution of the interface and the reason was the ease extension of the method from two space-dimensions to three space-dimensions. However,  we could also use other numerical representation techniques and in some cases an explicit representation of the interface may be beneficial due to a more accurate approximation of the interface.  

In future work, we aim at combining the method presented here with the method in \cite{HLZ16, Zah18} to also allow for surfactants.

\section*{Acknowledgement}
This research was supported by the Swedish Research Council Grant No. 2014-4804.
\section*{References}

%\bibliographystyle{elsarticle-num}
%%\biboptions{sort&compress}
%\bibliography{ref}

\begin{thebibliography}{10}
\expandafter\ifx\csname url\endcsname\relax
  \def\url#1{\texttt{#1}}\fi
\expandafter\ifx\csname urlprefix\endcsname\relax\def\urlprefix{URL }\fi
\expandafter\ifx\csname href\endcsname\relax
  \def\href#1#2{#2} \def\path#1{#1}\fi

\bibitem{UnTr92}
S.~O. Unverdi, G.~Tryggvason, A front-tracking method for viscous,
  incompressible, multi-fluid flows, J. Comput. Phys. 99~(1) (1992) 180--180.

\bibitem{SuSmOs94}
M.~Sussman, P.~Smereka, S.~Osher, A level set approach for computing solutions
  to incompressible two-phase flow, J. Comput. Phys. 114~(1) (1994) 146 -- 159.

\bibitem{FrBe10}
T.-P. Fries, T.~Belytschko, The extended/generalized finite element method: An
  overview of the method and its applications, Int. J. Numer. Meth. Engng
  84~(3) (2010) 253--304.

\bibitem{GrRe11}
S.~Gross, A.~Reusken, Numerical Methods for Two-phase Incompressible Flows,
  Springer Series in Computational Mathematics, Vol 40, 2011.

\bibitem{GaTo12}
S.~Ganesan, L.~Tobiska, Arbitrary {L}agrangian--{E}ulerian finite-element
  method for computation of two-phase flows with soluble surfactants, J.
  Comput. Phys. 231~(9) (2012) 3685 -- 3702.

\bibitem{BGN15}
J.~W. Barrett, H.~Garcke, R.~N{\"u}rnberg, A stable parametric finite element
  discretization of two-phase {N}avier--{S}tokes flow, J Sci Comput 63~(1) (2015)
  78--117.

\bibitem{GiHyFe18}
F.~Gibou, D.~Hyde, R.~Fedkiw, Sharp interface approaches and deep learning
  techniques for multiphase flows, J. Comput. Phys., in press (2018).

\bibitem{Win07}
C.~Winkelmann, Interior penalty finite element approximation of {N}avier--{S}tokes
  equations and application to free surface flows, Ecole Polytechnique Federale
  de Lausanne, PhD thesis (these no 3971), 2007.

\bibitem{GrRe07}
S.~Gross, A.~Reusken, An extended pressure finite element space for two-phase
  incompressible flows with surface tension, J. Comput. Phys. 224 (2007) 40 --
  58.

\bibitem{ScRaGrWa15}
B.~Schott, U.~Rasthofer, V.~Gravemeier, W.~A. Wall, A face-oriented stabilized
  {Nitsche}-type extended variational multiscale method for incompressible
  two-phase flow, Int. J. Numer. Meth. Engng 104~(7) (2015) 721--748.

\bibitem{KrBoSiVo}
R.~Kramer, P.~Bochev, C.~Siefert, T.~Voth, An extended finite element method
  with algebraic constraints {(XFEM-AC)} for problems with weak discontinuities,
  Comput. Methods Appl. Mech. Engrg. 266 (2013) 70 -- 80.

\bibitem{HaHa02}
A.~Hansbo, P.~Hansbo, An unfitted finite element method, based on {Nitsche's}
  method, for elliptic interface problems, Comput. Methods Appl. Mech. Engrg.
  191 (2002) 5537--5552.

\bibitem{BBH09}
R.~Becker, E.~Burman, P.~Hansbo, A {Nitsche} extended finite element method for
  incompressible elasticity with discontinuous modulus of elasticity, Comput.
  Methods Appl. Mech. Engrg. 198 (2009) 3352 -- 3360.

\bibitem{HaLaZa14}
P.~Hansbo, M.~G. Larson, S.~Zahedi, A cut finite element method for a {S}tokes
  interface problem, Appl. Numer. Math. 85 (2014) 90--114.

\bibitem{BuClHaLaMa15}
E.~Burman, S.~Claus, P.~Hansbo, M.~G. Larson, A.~Massing, Cut{FEM}:
  Discretizing geometry and partial differential equations, Int. J. Numer.
  Meth. Engng 104~(7) (2015) 472--501.

\bibitem{Bur10}
E.~Burman, Ghost penalty, C. R. Acad. Sci. Paris, Ser. I 348~(21-22) (2010)
  1217 -- 1220.

\bibitem{BH12}
E.~Burman, P.~Hansbo, Fictitious domain finite element methods using cut
  elements: {II. A} stabilized {Nitsche} method, Appl. Numer. Math. 62~(4)
  (2012) 328 -- 341.

\bibitem{WZKB}
E.~Wadbro, S.~Zahedi, G.~Kreiss, M.~Berggren, A uniformly well-conditioned,
  unfitted {Nitsche} method for interface problems, BIT Numer. Math. 53 (2013)
  791--820.

\bibitem{HanLarZah15}
P.~Hansbo, M.~G. Larson, S.~Zahedi, Stabilized finite element approximation of
  the mean curvature vector on closed surfaces, SIAM J. Numer. Anal. 53~(4)
  (2015) 1806--1832.

\bibitem{HLZ16}
P.~Hansbo, M.~Larson, S.~Zahedi, A cut finite element method for coupled
  bulk-surface problems on time-dependent domains, Comput. Methods Appl. Mech.
  Engrg. 307 (2016) 96 -- 116.

\bibitem{Zah18}
S.~Zahedi, A space-time cut finite element method with quadrature in time, in:
  Geometrically Unfitted Finite Element Methods and Applications, Lecture Notes
  in Computational Science and Engineering, Springer, 2018, pp. 281--306.

\bibitem{Le16}
C.~Lehrenfeld, High order unfitted finite element methods on level set domains
  using isoparametric mappings, Comput. Methods Appl. Mech. Engrg. 300 (2016)
  716 -- 733.

\bibitem{Nit}
J.~Nitsche, {Uber ein Variationsprinzip zur L{\"o}sung von Dirichlet-Problemen
  bei Verwendung von Teilr\"{a}umen, die keinen Randbedingungen unterworfen
  sind.}, Abh. Math. Sem. Univ. Hamburg 36 (1971) 9 -- 15.

\bibitem{BoiBur16}
T.~Boiveau, E.~Burman, A penalty-free {N}itsche method for the weak imposition
  of boundary conditions in compressible and incompressible elasticity, IMA
  Journal of Numerical Analysis 36~(2) (2016) 770--795.

\bibitem{GuOl18}
J.~Guzmán, M.~Olshanskii, Inf-sup stability of geometrically unfitted {S}tokes
  finite elements, Math. Comp 87 (2018) 2091--2112.

\bibitem{HanLarZah16}
P.~Hansbo, M.~G. Larson, S.~Zahedi, A cut finite element method for coupled
  bulk-surface problems on time-dependent domains, Comput. Methods Appl. Mech.
  Engrg. 307 (2016) 96--116.

\bibitem{CL15}
C.~Lehrenfeld, The {N}itsche {XFEM-DG} space-time method and its implementation
  in three space dimensions, SIAM J. Sci. Comput. 37~(1) (2015) A245 -- A270.

\bibitem{BuZu11}
E.~Burman, P.~Zunino, Numerical approximation of large contrast problems with
  the unfitted {Nitsche} method, in: Blowey J., Jensen M. (eds) Frontiers in
  Numerical Analysis - Durham 2010, Lecture Notes in Computational Science and
  Engineering, vol 85, Springer, Berlin, Heidelberg, 2011.

\bibitem{AnnHauDol12}
C.~Annavarapu, M.~Hautefeuille, J.~E. Dolbow, A robust {N}itsche’s
  formulation for interface problems, Comput. Methods Appl. Mech. Engrg.
  225-228 (2012) 44 -- 54.

\bibitem{OshFed01}
S.~Osher, R.~P. Fedkiw, Level set methods: An overview and some recent results,
  J. Comput. Phys. 169~(2) (2001) 463 -- 502.

\bibitem{Set01}
J.~Sethian, Evolution, implementation, and application of level set and fast
  marching methods for advancing fronts, J. Comput. Phys. 169~(2) (2001) 503 --
  555.

\bibitem{TryBuEsetal01}
G.~Tryggvason, B.~Bunner, A.~Esmaeeli, D.~Juric, N.~Al-Rawahi, W.~Tauber,
  J.~Han, S.~Nas, Y.-J. Jan, A front-tracking method for the computations of
  multiphase flow, J. Comput. Phys. 169~(2) (2001) 708 -- 759.

\bibitem{SuFa99}
M.~Sussman, E.~Fatemi, An efficient, interface-preserving level set
  redistancing algorithm and its application to interfacial incompressible
  fluid flow, SIAM J. Sci. Comput. 20~(4) (1999) 1165--1191.

\bibitem{MiGi07}
C.~Min, F.~Gibou, Geometric integration over irregular domains with application
  to level-set methods, J. Comput. Phys. 226 (2007) 1432--1443.

\bibitem{Dz88}
G.~Dziuk, Finite elements for the {B}eltrami operator on arbitrary surfaces,
  in: Partial differential equations and calculus of variations, Vol. 1357 of
  Lecture Notes in Math., Springer, Berlin, 1988, pp. 142--155.

\bibitem{Ba01}
E.~B{\"a}nsch, Finite element discretization of the {N}avier--{S}tokes equations
  with a free capillary surface, Numer. Math. 88~(2) (2001) 203--235.

\bibitem{Hys06}
S.~Hysing, A new implicit surface tension implementation for interfacial flows,
  Int. J. Numer. Meth. Fluids 51~(6) (2006) 659--672.

\bibitem{LarZah17}
M.~G. Larson, S.~Zahedi, Stabilization of high order cut finite element methods
  on surfaces, $ $\href {http://arxiv.org/abs/1710.03343}
  {\path{arXiv:1710.03343}}.

\bibitem{MuKuOb13}
B.~Müller, F.~Kummer, M.~Oberlack, Highly accurate surface and volume
  integration on implicit domains by means of moment-fitting, Int. J. Numer.
  Meth. Engng 96~(8) (2013) 512--528.

\bibitem{Sa15}
R.~Saye, High-order quadrature methods for implicitly defined surfaces and
  volumes in hyperrectangles, SIAM J. Sci. Comput. 37~(2) (2015) A993--A1019.

\bibitem{FriOmer16}
T.-P. Fries, S.~Omerović, Higher-order accurate integration of implicit
  geometries, Int. J. Numer. Meth. Engng 106~(5) (2016) 323--371.

\bibitem{OlRe18}
M.~A. Olshanskii, A.~Reusken, Trace finite element methods for PDEs on
  surfaces, in: Geometrically Unfitted Finite Element Methods and Applications,
  Lecture Notes in Computational Science and Engineering, Springer, 2018, pp.
  211--258.

\bibitem{example-Nasa}
S.~Tanveer, G.~Vascondelos, Time-evolving bubbles in two-dimensional {S}tokes
  flow, Journal of Fluid Mechanics 301 (1995) 325--344.

\bibitem{strain-flow-exact-sol}
S.~Richardson, Two-dimensional bubbles in viscous flow, Journal of Fluid
  Mechanics 33 (1968) 475--493.

\bibitem{Benchmark}
S.~Hysing, S.~Turek, D.~Kuzmin, N.~Parolini, E.~Burman, S.~Ganesan, L.~Tobiska,
  Quantitative benchmark computations of two-dimensional bubble dynamics, Int.
  J. Numer. Meth. Fluids 60~(11)  1259--1288.

\bibitem{Ran04}
R.~Rannacher, Incompressible Viscous Flows, American Cancer Society, 2004,
  Ch.~6.

\bibitem{Grace}
J.~Grace, Shapes and velocities of bubbles rising in infinite liquid,
  Transactions of the Institution of Chemical Engineers 51 (1973) 116--120.

\bibitem{Benchmark-3D}
Z.~Ge, J.-C. Loiseau, O.~Tammisola, L.~Brandt, An efficient mass-preserving
  interface-correction level set/ghost fluid method for droplet suspensions
  under depletion forces, J. Comput. Phys. 353 (2018) 435 -- 459.

\bibitem{Benchmark-3D_2}
M.~van Sint~Annaland, N.~G. Deen, J.~A.~M. Kuipers, Numerical simulation of gas
  bubbles behaviour using a three-dimensional volume of fluid method, Chemical
  Engineering Science 60~(11) (2005) 2999 -- 3011.

\end{thebibliography}
%\end{document}

\end{document}